\documentclass{article}
\usepackage{latexsym}
\usepackage{amsmath}
\usepackage{amsfonts}
\usepackage{epsfig}
\usepackage{stmaryrd}
\usepackage{amssymb}
\usepackage{pxfonts}
\usepackage{wasysym}
\usepackage{bm}

\newtheorem{theorem}{Theorem}

\newtheorem{lemma}{Lemma}

\usepackage{algorithmic}
\usepackage{algorithm}

\usepackage[colorlinks=true]{hyperref} 
\usepackage{fancyhdr}
\usepackage[normalem]{ulem}

\author{Francisco Bernal\footnote{CMAP, Ecole Polytechnique, email: {\tt francisco.bernal@polytechnique.edu}}}
\title{An implementation of Milstein's method for general bounded diffusions}

\date{}

\begin{document}
\maketitle

\begin{abstract}
Despite its generality and powerful convergence properties, Milstein's method for functionals of spatially bounded stochastic differential equations is widely regarded as difficult to implement. This has likely prevented it from being utilised in applications. In this paper, we design and analyse in detail one such implementation. The presented method turns out to be on par with other, popular schemes in terms of computational cost---but with a (nearly) linear weak convergence rate under the usual smoothness requirements on coefficients and boundary. Two byproducts of theoretical interest are a new, non-standard rank-one update formula, and a connection between numerics of bounded diffusions and Eikonal equations. Three examples are worked out, confirming the accuracy and robustness of the method.
  
\end{abstract}

(Accepted for publication in the Journal of Scientific Computing, 2018)\newline



\section{Introduction}\label{S:Introduction}
\paragraph*{Scope and motivation.} We are interested in {\em efficient} weak schemes for Feynman-Kac functionals of stochastic differential equations (SDEs) related to {\em general} linear boundary value problems (BVPs). By efficiency, we refer to the computational cost required to bound the expected error within a given tolerance and a given confidence interval. (This depends on the weak rate of convergence w.r.t. the timestep $h$ of the SDE scheme, and on its cost per iteration.) By general, we mean that no restrictions are posed on the coefficients and boundary beyond those required for solvability of the associated SDE (i.e. a minimum degree of regularity of both). In particular, we cover BVPs with mixed boundary conditions (BCs) (i.e. of Dirichlet type on some portion of the boundary and Robin elsewhere), for which the stochastic representation is in terms of functionals of SDEs spatially enclosed by reflecting and stopping boundaries.  

Such stochastic schemes are the basis for the pointwise solution, $u({\bf x}_0,t_0)$, of linear BVPs in ${\mathbb R}^D$ via the Monte Carlo method \cite{Milstein_Tretyakov_Book}. Besides the mathematical interest of this connection by itself, the Monte Carlo approach is computationally advantageous in many applications, v.g.   \cite{IterPDD,Bernal2014,Cao2007,Maire&Simon,Mascagni2002}. 
The stochastic numerics presented here can also be adapted to construct a variety of probability densities of a Wiener diffusion inside a domain enclosed by absorbing and/or reflecting surfaces, including the survival probability, the time spent close to the reflecting boundary, and first-passage times from the domain  \cite{Licata&Grill,Redner2001,Schwabedal_PRL}.

\paragraph*{State of the art.} Given a general linear BVP, it is challenging to numerically integrate its associated Feynman-Kac functionals with a satisfactory weak order of convergence w.r.t. $h$ ($\delta$ henceforth). This is specially true when the BVP has mixed BCs---or equivalently,  when the SDE is stopped on some portion of the boundary and reflected elsewhere. The method of Constantini {\em et al.} \cite{Constantini1998} may be considered as the standard; it is very easy to code but has a (proven) $\delta$ of just up to one half. The difficulties pertain to the determination (in a weak sense) of the first-exit point (in the presence of stopping boundaries---see \cite{OurReview} and references therein), and of the local time (in the presence of reflecting ones). 
Other schemes have a better $\delta$, but they are less general; the following list is not exhaustive. Gobet's half-space approximation \cite{Gobet2001_BB} is  very useful in practice. He has theoretically analised the schemes for either killed (rather than stopped) {\em or} reflected diffusions; and proved that many Feynman-Kac functionals converge linearly with $h$---including those for homogeneous Dirichlet BCs, and those for Neumann (rather than Robin) BCs. Bossy {\em et al.} \cite{BossyGobetTalay} put forward a straightforward method which has $\delta=1$ for the Feynman-Kac functionals associated with homogeneous Neumann BCs. A relatively new approach (see \cite{Bayer2010} and references therein) relies on adaptive $h$, based on a running error estimate. It has shown very good results (including in nonconvex domains) when the diffusion is a Brownian motion, but we are not aware of extensions to more general cases. Also sticking to the Brownian motion, the method of randomisation can deal accurately with all sorts of linear BCs (see \cite{Maire&Simon} and references therein).

Besides the method of Constantini {\em et al.}, there is another one, due to Milstein, which also is completely general {\em and} has a proven ${\cal O}(h)$ weak convergence rate \cite{Milstein_Neumann,Milstein_WoE}: we will call it MM in short. Paradoxically (since it predates all the methods quoted above), MM has gone underreported in the literature, most likely because of a reputation of being complex to implement. (In fact, we are not aware of numerical experiments with it having been published, beyond one example in \cite[p. 372]{Milstein_Tretyakov_Book} and a few more in our review \cite{OurReview}, both on purely stopped diffusions.) For instance, MM is described as follows in \cite[p. 77]{Constantini1998}: "(...) two other weak discretization schemes (...) are considered in domains with smooth boundary (...) One of them achieves the rate of convergence $h$ but, as pointed out by the author, is difficult to implement". Moreover, quoting from \cite[p. 280]{Gobet2001_BB}: "An appropriate Markov chain approximation at random discretization times has been studied by Milstein (...) His procedure requires at each step near the boundary to change coordinates, and, by the way, the algorithm seems to be difficult to implement".

\paragraph*{Our contribution.} The reason why MM is hard to implement is that the algorithms are incomplete. At a given point, the user is instructed to "take the ellipsoid tangent to the boundary", or to perform a given rotation or a certain matrix decomposition close to reflecting boundaries; but details as to how to are glossed over both in the seminal papers as in the later book \cite[chapter 6]{Milstein_Tretyakov_Book}. In this paper, we provide a practical and theoretically sound such implementation of MM---apparently the first one.
(Henceforth, we will call it ``our implementation'' or "Algorithm \ref{A:WoE}'', to distinguish it from 
the original algorithm.\footnote{The Matlab code and data files used for this paper are available at the journal repository, and upon request from the author.})

Regarding the determination of the tangent ellipsoid (Section \ref{SS:Ellipsoid}), we have established a novel link with an anisotropic Eikonal equation (Lemma \ref{Th:Anisotropic}) which allows one, theoretically, to do it exactly. Nonetheless, this approach is unpractical except with a constant diffusion matrix (as in Example \ref{SS:Example3}). Therefore, a fast ${\cal O}(D^2)$ half-space approximation (Lemma \ref{Th:Tangent_ellipsoid}) is also provided. For the reflections (Section \ref{SS:Rotation}), we have introduced a non-standard, non-Cholesky update formula which endowes the overall implementation with ${\cal O}(D^2)$ complexity (Algorithm \ref{A:Update}). This is relevant because MM is designed for Monte Carlo simulations, which are often meant for high-dimensional problems ($D\gg 1$). In addition to those three main contributions, we discuss every other aspect and provide fast 
recipes for them. In sum, this paper lifts what were the main obstacles to the routinary utilization of Milstein's method for general bounded diffusions.

A preliminary, simplified version of the implementation presented here, valid only for smooth, purely stopping boundaries, was sketched in the review paper \cite{OurReview}. That algorithm, however, is now superseded by Algorithm \ref{A:WoE}. 

Our implementation turns out to have a cost per time step comparable (less than $2$ times larger) to that of the integrator by  Constantini {\em et al.}---however, with a (nearly) linear weak convergence rate. To put this fact in perspective, the complexity of a Monte Carlo simulation (proportional to the CPU time) with error tolerance $0<\epsilon\ll1$ drops dramatically (by a factor $1/\epsilon$) from $\delta=1/2$ to $\delta=1$ \cite{Giles_y_yo}. Substantial further gains would be possible by combining Algorithm \ref{A:WoE} with Multilevel \cite{Giles_y_yo}, extrapolation \cite{Sara}, or both \cite{Pages}. Critically, all of the previous strategies rely on the {\em a priori} knowledge of $\delta$---thus highlighting the suitability of MM thanks to its sound theoretical foundation.  

While MM has an asymptotically proven $\delta=1$ under the assumptions of adequate smoothness of the coefficients and of the boundary, our implementation might not exactly reproduce it in two cases:
\begin{enumerate}
	\item Unless $\Omega$ has a trivial shape (or is a combination thereof), the distance map will itself be a numerical approximation---thus, possibly spoiling $\delta=1$. (Note, however, that the same holds for {\em any} numerical scheme for bounded diffusions, not just the one presented in this paper.)
\item When the absorbing boundaries are curved and the fast recipe in Lemma \ref{Th:Tangent_ellipsoid} for constructing the smallest tangent ellipsoid is employed (as will be typically the case). The reason is that Lemma \ref{Th:Tangent_ellipsoid} relies on locally approximating the boundary by its tangent plane, but we have not rigourously proved convergence of $\delta\to 1$ as $h\to 0$. 
\end{enumerate}  
Consequently, we claim that our algorithm has a "(nearly) linear weak convergence rate". On the other hand, we emphasize that all of the numerical experiments involving smooth boundaries which we have carried out univocally suggest $\delta=1$, in practice.

We close the Introduction by briefly commenting on preprocessing.

The main effort which must be independently undertaken before the Monte Carlo simulation is the generation of a signed distance map for $\Omega$. Again, we stress that this is required for any numerical scheme for bounded SDEs---although it often goes unmentioned. In the simplest cases, like a ball or a parallellepiped, an exact distance formula is available. Otherwise, we propose solving an Eikonal equation via the Fast Marching Method (check Section \ref{SS:Eikonal}).

Finally, if the exact tangent ellipsoid is needed, the anisotropic Eikonal equation in Lemma \ref{Th:Tangent_ellipsoid} must be numerically solved in advance. Example \ref{SS:Example3}  has been crafted to illustrate a situation where this would be highly advantageous, for that problem is very hard to solve with other numerical schemes.\newline
Since preprocessing is either not particular to the algorithm presented here, or up to some point optional---and in both cases to be tackled with an independent method---its cost has not been explicitly included in the main discussion.

\paragraph*{Organisation of the paper.} Section \ref{S:Representation} recalls the theoretical connection between linear second-order BVPs with mixed BCs, and bounded stopped/reflected SDEs. (The form of the Feynman-Kac formulas there is not the most usual one, but one tailored to MM.) In Section \ref{S:Algorithm}, the original MM is described, and the new Algorithm \ref{A:WoE} is listed. Implementation details are discussed around several new lemmas in Section \ref{S:Implementation}, which is the core of the paper. Three numerical examples are worked out in Section \ref{S:Experiments}; and conclusions are drawn in Section \ref{S:Conclusions}. To avoid clutter, all proofs of the lemmas in Section \ref{S:Implementation} have been moved into Appendix \ref{Ap:Proofs}. Finally, Appendix \ref{Ap:Code} lists a few relevant Matlab code snippets.

\section{Feynman-Kac formulas in Milstein's form}\label{S:Representation}

Let $D\geq 2$, $\Omega\subset {\mathbb R}^D$ be a bounded domain, and $\Omega={\overline{\Omega}\cup\partial\Omega}$, where the open connected set ${\overline\Omega}$ is the interior of the domain and $\partial\Omega$ its boundary. Consider the linear parabolic BVP of second order with mixed BCs:

\begin{equation}
\label{eq:ParabolicBVP_mixedBCs}
\left\{
\begin{array}{ll}
\frac{\partial u}{\partial t}= {\cal L}({\bf x},t) u + c({\bf x},t)u + f({\bf x},t) & \textrm{ if } 0<t\leq T, {\bf x}\in{\overline\Omega},\\
u= p({\bf x}) & \textrm{ if } t=0, {\bf x}\in\Omega,\\
u= g({\bf x},t) & \textrm{ if } 0<t\leq T, {\bf x}\in\partial\Omega_A,\\
\frac{\partial u}{\partial N}= {\varphi}({\bf x},t)u + \psi({\bf x},t) & \textrm{ if } 0<t\leq T, {\bf x}\in\partial\Omega_R,
\end{array}
\right.
\end{equation}

where $T>0$, $\varphi({\bf x},t)\leq 0$, 
and the differential generator is given by

\begin{equation}
{\cal L}({\bf x},t)u= \frac{1}{2}\sum_{i,j=1}^D a_{ij}({\bf x},t)\frac{\partial^2 u}{\partial x_i\partial x_j} + \sum_{k=1}^D b_k({\bf x},t)\frac{\partial u}{\partial x_k}.
\end{equation}

The matrix $A({\bf x},t):=[a_{ij}]_{i,j=1}^D$ is positive definite, and ${\bf b}({\bf x},t):=(b_1,\ldots,b_D)^T$ is called the drift. All of the coefficient functions in (\ref{eq:ParabolicBVP_mixedBCs}), namely $a_{ij},b_i,c,f,p,g,\varphi$ and $\psi$ are assumed continuous, and complying with the compatibility conditions at time $t=0$ (see \cite{Rusos} or \cite[equations (2.15)-(2.17)]{Constantini1998}). The boundary is decomposed as $\partial\Omega=\partial\Omega_A\cup\partial\Omega_R\cup\partial\Omega_S$, such that $\partial\Omega_A\cap\partial\Omega_R=\partial\Omega_A\cap\partial\Omega_S=\partial\Omega_R\cap\partial\Omega_S=\emptyset$.
The outward\footnote{In the SDE literature, the normal is usually taken inwards. Here we follow the PDE convention.} unit normal vector ${\bf N}$ is assumed to be well defined on the boundary save perhaps on a set $\partial\Omega_S$; 
$\partial\Omega_A$ stands for the portion of the boundary (if any) where Dirichlet BCs are imposed; and on $\partial\Omega_R$, BCs involving the normal derivative, (i.e.  ${\bf N}^T\nabla u$) hold, where $\nabla u=(\partial u/\partial x_1,\ldots,\partial u/\partial x_D)^T$. (BCs involving oblique derivatives will not be considered in this paper.) Such BCs are of Neumann type iff $\varphi=0$, or of Robin type otherwise.

\paragraph*{Sufficient conditions for existence of a unique classical solution to (\ref{eq:ParabolicBVP_mixedBCs}).}
When $\partial\Omega_A=\emptyset$ (respectively $\partial\Omega_R=\emptyset$) we say the BCs are purely reflecting (resp. purely stopping), while 
when both $\partial\Omega_A \neq \emptyset \neq \partial\Omega_R$, we say the BCs are mixed. If the BCs are purely reflecting (resp. purely stopping), and $\partial\Omega_S=\emptyset$,  theorem 2.6 (resp. theorem 2.7) in \cite{Constantini1998} (see also \cite{Rusos}) ensure the existence of a unique classical solution, with regularity dependent on that of the BVP coefficients and of $\partial\Omega$. (By a classical solution, we mean that $u({\bf x},t)$ lives in the H{\"o}lder space ${\cal C}^{1,2}([0,T]\times{\overline{\Omega}})$.) In the mixed BCs case, this connection is less general and more dependent on the smoothness of the boundary \cite{Lieberman1986,Miranda1970}.

\paragraph*{Stochastic representation of the pointwise solution to (\ref{eq:ParabolicBVP_mixedBCs}).} Under slightly stronger conditions, the stochastic representation of the PDE with mixed BCs (\ref{eq:ParabolicBVP_mixedBCs}) expresses its pointwise solution $u({\bf x}_0,t)$ as the expected value of a  functional of an SDE starting at ${\bf x}_0$ at time zero and being normally reflected on $\partial\Omega_R$ and stopped on $\partial\Omega_A$. 

Let $\sigma({\bf x},t)$ (called the diffusion matrix) be defined by $\sigma({\bf x},t)\sigma^T({\bf x},t)=A({\bf x},t)$. (This is always possible since $A$ is positive definite. Hence, $\det{\sigma}\neq 0$.) The following result is an extension of the well-known Feynman-Kac formulas for purely reflected ($\partial\Omega_A=\emptyset$) and purely stopped ($\partial\Omega_R=\emptyset$) diffusions, adapted from \cite[theorem 2.5]{Constantini1998} and \cite[chapter 6]{Milstein_Tretyakov_Book}.

\begin{theorem}
	Assume that: i) a classical unique solution to (\ref{eq:ParabolicBVP_mixedBCs}) does exist; ii) there exists a constant $L>0$ such that for ${\bf x},{\bf y}\in{\overline\Omega}$ and $t\in[0,T]$
	\begin{eqnarray}
	||\sigma({\bf x},t)-\sigma({\bf y},t)||\leq L||{\bf x}-{\bf y}||, \\
	||{\bf b}({\bf x},t)-{\bf b}({\bf y},t)||\leq L||{\bf x}-{\bf y}||;
	\end{eqnarray} 
	iii) $\partial\Omega$ is piecewise ${\cal C}^1$ (i.e. ${\cal C}^1$ save on maybe a set $\partial\Omega_S$); and iv) either  $\overline{\Omega}$ is convex, or $u({\bf x},t)$ can be extended to a function ${\cal C}^{1,2}([0,T]\times{\mathbb R}^d)$. Then, for $0\leq t\leq T$ the pointwise solution of (\ref{eq:ParabolicBVP_mixedBCs}) admits the following stochastic representation:
	
	\begin{equation}
	\label{eq:Parabolic_representation}
	u({\bf x}_0,t)= {\mathbb E}[\phi_{\tau}]:= {\mathbb E}\big[\, q({\bf X}_{\tau})Y_{\tau}+Z_{\tau} \,\big], 
	\end{equation}
	
	where
	
	\begin{equation}
\label{eq:Milstein_IC_and_BC}	
	q({\bf X}_{\tau})=\left\{
	\begin{array}{ll}
	g({\bf X}_{\tau},T-\tau), & \textrm{ if } \tau<T,\\
	p({\bf X}_T), & \textrm{ if } \tau\geq T, 
	\end{array}
	\right.
	\end{equation}
	and the processes $({\bf X}_t,Y_t,Z_t,\xi_t)$ are governed by the following set of SDEs driven by a standard $D-$dimensional Wiener process ${\bf W}_t$:  
\begin{equation}
\label{eq:Milstein_sys}
\left\{
\begin{array}{ll}
d{\bf X}_t= \Big( {\bf b}({\bf X}_t,T-t) - \sigma({\bf X}_t,T-t){\bm\mu}\Big)dt + \sigma({\bf X}_t,T-t)d{\bf W}_t -{\bf N}({\bf X}_t)d\xi_t & {\bf X}_0={\bf x}_0,\\
dY_t= c({\bf X}_t,T-t)Y_t dt + {\varphi}({\bf X}_t,T-t) Y_t d\xi_t + Y_t{\bm\mu}^Td{\bf W}_t& Y_0=1,\\ 
dZ_t= f({\bf X}_t,T-t)Y_t dt + \psi({\bf X}_t,T-t) Y_t d\xi_t + Y_t{\bf F}^Td{\bf W}_t & Z_0= 0,\\
d\xi_t= {\bf 1}_{\{{\bf X}_t\in\partial\Omega_R\}}dt & \xi_0=0.
\end{array}
\right.
\end{equation}
	
	Above, ${\bf 1}_{\{H\}}$ is the indicator function ($1$ if $H$ is true and $0$ otherwise); $\tau=\inf_{t}\{{\bf X}_t\in\partial\Omega_A\}$ is the first exit time  (or first passage time) from $\Omega$; which takes place at the first exit point ${\bf X}_{\tau}\in\partial\Omega_A$; and $\xi_t$ is called the local time.
\end{theorem}
The functions ${\bm \mu}:\Omega\mapsto{\mathbb R}^D$  and ${\bf F}:\Omega\mapsto{\mathbb R}^D$ ---the former a consequence of Girsanov's theorem and the latter of the expectation of Ito's integral being zero---
are at this point rather arbitrary.
However, properly choosing the function ${\bm\mu}$ will later be crucial for the WoE numerical scheme.

In the remainder of the paper, ${\bf F}$ is set to zero. (See \cite{IterPDD} for an application where it is not.) 

Many expected values pertaining to a population of particles can be accommodated into this framework---see \cite{Constantini1998}. For instance, if $p=g=c=\varphi=\psi=0$, $f=1$, then $u({\bf x}_0,t)$ is the mean absorption time of a particle starting at ${\bf x}_0$ and driven by the drift ${\bf b}$ and diffusion matrix $\sigma$ in the interval $[0,T]$. Analogously, if $p=g=c=\varphi=f=0$, $\psi=1$, $u({\bf x}_0,t)$ is the mean local time. 

\paragraph*{Elliptic equations.} Equation (\ref{eq:ParabolicBVP_mixedBCs}) can be formally 
transformed into an elliptic BVP with mixed BCs by assuming that $\partial u/\partial t=0$, thus dropping the dependence on time from $u$ and all the coefficients; letting $T\shortrightarrow \infty$; and dropping the initial condition $p$. If $\partial\Omega_R=\emptyset$ (purely stopped diffusions / Dirichlet BCs), the stochastic representation derived from Theorem 1 still holds as long as $c\leq 0$ and ${\mathbb E}[\tau]<\infty$ \cite[section 4.4.5]{Gobet2016}. (Note that in the purely reflected case, the latter condition would be impossible.) To the best of our knowledge, there is no rigourously proved stochastic representation for elliptic BVPs with mixed BCs.  Therefore, we will assume in the remainder of this paper that, if: i) the time-independent equivalent conditions of those in Theorem 1 are in place; ii) $c\leq 0$; and iii) ${\mathbb E}[\tau]<\infty$; then the time-independent equivalent representation given by Theorem 1 holds.

\section[Overview of WoE]{Overview of Milstein's method and the proposed implementation (Algorithm 1)}\label{S:Algorithm}

\subsection{Notation}

In order to approximate the SDE system (\ref{eq:Milstein_sys}) numerically, a small, constant timestep $h>0$ is set, and each random realization of $({\bf X}_t)_{0\leq t\leq \min{(\tau,T)}}$ is replaced by a chain (which we may call a ``trajectory'') ${\bf\hat X}_0={\bf x}_0,{\bf\hat X}_1,\ldots,{\bf\hat X}_K$, where ${\bf\hat X}_k\approx {\bf X}_{kh}$, $k=0,1,...,K$, ${\hat\tau}:=Kh\approx \tau$, and $t_k:=kh$. The chains $\{{\hat Y}_k\}_{k=0}^K$, $\{{\hat Z}_k\}_{k=0}^K$ and $\{{\hat \xi}_k\}_{k=0}^K$ are defined analogously. In order to lighten the notation, we drop the hats and ${\bf X}_k$, $Y_k$ etc. are to be understood as the discretized counterparts ${\hat{\bf X}}_k$, ${\hat Y}_k$, etc. unless stated otherwise. In general, functions with subindex $k$ are evaluated at $({\bf X}_k,T-t_k)$, like v.g. $c_k=c({\bf X}_k,T-t_k)$.

Vectors are by default column vectors, and $||\cdot||$ is always the Euclidean norm. For ${\bf x}\in{\mathbb R}^D$ and $\partial\Omega$ smooth, let us define: 
\begin{itemize}
\item ${\bf x}^{\partial\Omega}=\arg\min_{{\bf y}\in\partial\Omega}||{\bf x}-{\bf y}||$ (i.e. the closest point on the boundary).
\item $|d({\bf x})|=||{\bf x}-{\bf x}^{\partial\Omega}||$ (i.e. the Euclidean distance to the boundary).
\item ${\bf N}({\bf x})={\bf N}({\bf x}^{\partial \Omega})$ (i.e. the normal vector pointing outwards).
\item $\Pi({\bf x})=\{{\bf y}\,|\,{\bf N}^T({\bf x}) ({\bf y}-{\bf x}^{\partial\Omega})=0\}$ (i.e. the hyperplane tangent to $\partial\Omega$ at ${\bf x}^{\partial\Omega}$).
\end{itemize}

We use the shorthand notation $|d_k|:=|d({\bf X}_k)|$, ${\bf N}_k:={\bf N}({\bf X}_k)$, and $\Pi_k:=\Pi({\bf X}_k)$.

$B({\bf x},R)$ stands for a $D-$dimensional ball of radius $R$ centred at ${\bf x}$, and $\partial B({\bf x},R)$ for its surface. Let $M$ be a $D\times D$ symmetric positive definite matrix with (real) spectrum given by $M{\bf m}_i=\mu_i{\bf m}_i$, $\mu_1\geq \mu_2\geq\ldots\geq\mu_D>0$. For $\rho>0$, define the  $D-$dimensional ellipsoid centred at ${\bf x}$ and oriented according to $M$ as
\begin{equation}
\label{eq:Elipsoide}
E^{\rho}_M({\bf x})= \{{\bf y}\,|\,({\bf y}-{\bf x})^TM^{-1}({\bf y}-{\bf x})=\rho^2\}.
\end{equation}

The major semiaxis of $E^{\rho}_M({\bf x})$ has length $\rho\mu_1^{1/2}$ and points along ${\bf m}_1$
\footnote{The eigenvalues and eigenvectors of $M^{-1}$ are the inverse and the same as those of $M$, respectively.}
, and so successively until the minor semiaxis, which runs along ${\bf m}_D$ with length $\rho\mu_D^{1/2}$. When $M=A({\bf X}_k,T-t_k)$, we may simply write $E^{\rho}_k:=E^{\rho}_{A({\bf X}_k,T-t_k)}({\bf X}_k)$. 

In particular, there is one value of $\rho$, which we call ${\overline{\rho}}$, such that $E^{\overline{\rho}}_M({\bf x})$ is tangent to $\partial\Omega$. We denote it as $\Sigma_M({\bf x}):=E^{\overline{\rho}}_M({\bf x})$---respectively, $\Sigma_k:=E^{\overline{\rho}}_k$. The surfaces of those ellipsoids are denoted as $\partial E^{\rho}_M({\bf x})$, $\partial E_k^{\rho}$, $\partial \Sigma_M({\bf x})$, and $\partial \Sigma_k$.

The eigenvalues of $A({\bf x},T-t)$ are $\lambda_1({\bf x},T-t)\geq \lambda_2({\bf x},T-t)\geq\ldots\geq \lambda_D({\bf x},T-t)>0$ (or simply $\lambda_1\geq\ldots\geq\lambda_D$ if the context permits.)

The notation $\omega\sim W$ means that $\omega$ is one realization drawn from a distribution $W$. In particular, ${\cal S}_D$ is the uniform distribution of points over $\partial B({\bf 0},1)$ (recall it is $D-$dimensional);  ${\cal B}$ is the distribution taking $\pm 1$ with equal probability; and ${\cal N}(m,s^2)$ is the Gaussian with mean $m$ and variance $s^2$.

\subsection{Description of Milstein's method}

\paragraph*{Remark.} For the sake of clarity, let us emphasize that what we refer to as MM is the combination of two different algorithms: one for purely stopped diffusions, called ``Walk on Ellipsoids'' (WoE) by Milstein, \cite{Milstein_WoE}\cite[section 6.3.2]{Milstein_Tretyakov_Book}; plus the reflection off $\Omega_R$ handled according to Milstein's one-step scheme \cite{Milstein_Neumann}\cite[section 6.6]{Milstein_Tretyakov_Book} for purely reflected diffusions. Under sufficient smoothness of the boundary and the coefficients (leading to a $C^4([0,T]\times{\overline\Omega})$ solution of \eqref{eq:ParabolicBVP_mixedBCs}), both algorithms enjoy proven ${\cal O}(h)$ weak convergence separately. The combination of both schemes (called MM by us) is theoretically analised (and nowhere tested) in \cite[theorem 6.2]{Milstein_Tretyakov_Book}; under the same contraints as above, it enjoys linear weak convergence rate as well. Hence, it furnishes a method for stopped/reflected SDEs and hence suitable for Monte Carlo solutions of BVPs with mixed BCs.  \newline

Let us now explain MM, starting with its first ingredient, WoE. In contrast with Euler-like integrators,  
WoE takes bounded steps in order to avoid overshooting the boundary. 
When ${\bf X}_k$ is not interacting with it (in a sense that will be clarified in a moment), the drift is first removed from (\ref{eq:Milstein_sys}) by setting ${\bm\mu}=\sigma^{-1}{\bf b}$. 
Then, WoE takes ${\bf X}_{k+1}={\bf X}_k+r\sigma{\bm\omega}$, where $r=\sqrt{Dh}$ and ${\omega\sim{\cal S}_D}$, whence 
\begin{equation}
\label{eq:Elipsoide}
||\bm\omega||^2 = 1 = \frac{1}{r^2} \big( {\bf X}_{k+1}-{\bf X}_k\big)^TA^{-1}_k\big( {\bf X}_{k+1}-{\bf X}_k\big). 
\end{equation}

By (\ref{eq:Elipsoide}), the possible values of ${\bf X}_{k+1}$ are distributed over $\partial E^r_k$.
As long as $t_k<T$, WoE has the sequence $\{{\bf X}_k\}_{k=1,2,...}$ hop over $\partial E^r_0$, $\partial E^r_1,\ldots$ until coming close enough to the boundary (at $k=i$, say) that ${\bf X}_i$ iteracts with it. 
Then, one of the following happens: i) a reflection on $\partial\Omega_R$ takes place, yielding ${\bf X}_{i+1}$; ii) a special ellipsoid must be taken for ${\bf X}_{i+1}$; or iii) $\partial\Omega_A$ is deemed hit and the trajectory stopped.  

If ${\bf X}_i$ is so close to the absorbing boundary that $E^r_i$ may intersect it (but still farther than $r^2=Dh$), WoE dictates that the tangent ellipsoid $\Sigma_i$ be taken instead for ${\bf X}_{i+1}$---in order to ensure that the trajectory cannot trespass the boundary. This may happen only when $|d_i|<r\sqrt{\lambda_1({\bf X}_i,T-t_i)}$. Eventually, either the trajectory is stopped at the projection on $\partial\Omega_A$ (this happens when the distance to it is less than $r^2$), or the initial condition is hit (i.e. $t_k\geq T$). 

If ${\bf X}_i$ is closer than $r$ to $\partial\Omega_R$, the one-step scheme in \cite[section 6.6]{Milstein_Tretyakov_Book} is used to handle the normal reflection on the boundary. This  involves a change of coordinates such that the new first component is pointing along $-{\bf N}_i$.

In a nutshell, Milstein's MM is structured in the following way:
\begin{enumerate}
	\item If $t_k\geq T$, the trajectory has 'reached' the initial condition without being stopped by $\partial\Omega_A$. Then, read the initial condition and finish.
	\item If not, and the trajectory is closer than $r^2$ to $\partial\Omega_A$, deem it stopped. Read the Dirichlet BC and finish.
	\item If neither of the above, and the trajectory is closer than $r$ to  $\partial\Omega_R$, perform the reflection, and continue.
	\item Otherwise, hop onto the surface of $E_k^r$ for ${\bf X}_{k+1}$---if that ellipsoid is fully contained in $\Omega$. If not, hop onto the surface of $\Sigma_k$, and continue.
\end{enumerate}

As it stands, MM cannot be used, for several points must be clarified first:
\begin{itemize}
	\item the (fast) determination of the boundary data: $|d({\bf x})|$, ${\bf N}({\bf x})$, ${\bf x}^{\partial \Omega}$, and $\Pi({\bf x})$.
	\item the (fast) determination of $\lambda_1({\bf x},T-t)$ for general diffusions,
	\item the (fast) determination of $\Sigma_{A({\bf x},T-t)}({\bf x})$ for general domains, and
	\item the rotation needed for handling the reflections (in a fast way).
\end{itemize} 
By "clarified'', we mean that specific recipes must be put in place for each of the above points. Very importantly, they must be cost-efficient, since they will be performed at potentially many time steps with each of the $N\gg 1$ realizations involved in the Monte Carlo estimate of the expectation (\ref{eq:Parabolic_representation}). 

We prove that all of the points above can be determined at ${\cal O}(D^2)$ cost without further assumptions than those for MM. Therefore, Algorithm \ref{A:WoE}---which is our implementation of MM and the main result of this paper---also has an overall complexity ${\cal O}(D^2)$ per time step, which is the minimum attainable order: that of computing the matrix-vector product 
$\sigma d{\bf W}_t$
\footnote{The product of a lower triangular matrix by a vector involves $1$ (first row) $+2$ (second row) $+\ldots+D$ (bottom row) $=D(D+1)/2$ multiplications and as many additions; i.e. an ${\cal O}(D^2)$ cost.} in \eqref{eq:Milstein_sys}. For that reason, we claim that our algorithm is "fast".

\begin{algorithm}
	\caption{A practical implementation of MM for general bounded diffusions}
	\begin{algorithmic}[1]
		\STATE{{\bf Data:} $1\gg h>0$, $r=\sqrt{Dh}$, a signed distance map $d({\bf x})$ such that $d({\bf x}\in{\overline\Omega})<0$}
		
		\STATE{Let ${\bf X}_0={\bf x}_0\in{\overline\Omega},Y_0=1,Z_0=0,\xi_0=0$, $t_0=0$, $d_0<0,{\bf X}_0^{\partial\Omega}, {\bf N}_0$, and $k=0$}
		
		\WHILE{neither the initial condition nor $\partial\Omega_A$ have been hit}
		
		\IF{$t_k\geq T$ (initial condition hit)}
		\STATE{Let ${\bf X}_T={\bf X}_k, Y_T=Y_k, Z_T=Z_k,\xi_T=\xi_k$, evaluate $p({\bf X}_T)$ in (\ref{eq:Milstein_IC_and_BC}) and {\bf finish}}
		
		\ELSIF{$|d_k|\leq r^2$ {\bf and} ${\bf X}_k^{\partial\Omega}\in\partial\Omega_A$ (absorption)}
		
		\STATE{$\,\,$Let $\tau=t_k$, ${\bf X}_{\tau}={\bf X}_k^{\partial\Omega}$, $Y_{\tau}=Y_k$, $Z_{\tau}=Z_k$, $\xi_{\tau}=\xi_k$, evaluate $g({\bf X}_{\tau},T-\tau)$ in (\ref{eq:Milstein_IC_and_BC}) and {\bf finish}}
		
		\ELSIF{$|d_k|\leq r$ {\bf and} ${\bf X}_k^{\partial\Omega}\in\partial\Omega_R$}
\STATE{Evaluate $c_k,f_k$ at $({\bf X}_k,T-t_k)$ ; let ${\bar \varphi}=\varphi({\bf X}_k^{\partial\Omega},T-t_k)$, ${\bar \psi}=\psi({\bf X}_k^{\partial\Omega},T-t_k)$}	
		\STATE{Compute Givens entries $\{\cos{\theta_2},\ldots,\sin{\theta_D}\}$ (Section \ref{SS:Rotation} and App. \ref{Ap:Code})}
		\STATE{Compute $\hat{{\bf b}}=Q{\bf b}({\bf X}_k,T-t_k)$,  $({\hat\varphi}'_1,\ldots,{\hat\varphi}'_D)^T=Q\nabla\varphi({\bf X}_k^{\partial\Omega},T-t_k)$ and $({\hat\psi}'_1,\ldots,{\hat\psi}'_D)^T=Q\nabla\psi({\bf X}_k^{\partial\Omega},T-t_k)$ by (\ref{F:Q_matrix_vector})}
		\STATE{Draw ${\vec \nu}=(\nu_1,\ldots,\nu_{D-1})$ with $\nu_i\sim {\cal B}$}
		
		\STATE{Compute $\Lambda_*{\vec\nu}$ according to Algorithm \ref{A:Update}, using $\sigma({\bf X}_k,T-t_k)$}
		
		\STATE{Compute $({\hat A}_{11},\ldots,{\hat A}_{D1} )^T$, the first column of $QA({\bf X}_k,T-t_k)Q^T$, by (\ref{F:A1})}
		
		\STATE{Let ${\bm\chiup}=(\chiup_1,\chiup_2,\ldots,\chiup_D)^T=
			\left\{\begin{array}{l}
			\chiup_1=\sqrt{{\hat A}_{11}r^2 + d_k^2}-|d_k|-{\hat b}_1r^2,\\
			\chiup_i=-{\bar\varphi}{\hat A}_{1i}r^2+(\Lambda_*{\vec \nu})_{i-1}r\textrm{, for } i=2,\ldots,D
			\end{array}\right.$}

		\STATE{Perform the reflection according to \cite[theorem 6.1]{Milstein_Tretyakov_Book}:
			\begin{eqnarray}
			\left\{\begin{array}{l}
			{\bf X}_{k+1}= {\bf X}_k + Q^T({\bm \chiup}+\hat{{\bf b}}r^2)\qquad\textrm{ (using (\ref{F:Q_matrix_vector}), see also Appendix \ref{Ap:Code})} \\
			
			Y_{k+1}= Y_k + \Bigg( c_kr^2 + 
			{\bar\varphi}\big(1+{\bar\varphi}|d_k|\big)\chiup_1 - \big(\sum_{j=2}^D{\hat A}_{1j}{\hat\varphi}'_j\big)r^2 + {\bar\varphi}^2\chiup_1^2 \Bigg)Y_k,\\
			
			Z_{k+1}= \,Z_k+\Bigg( f_kr^2 + {\bar\psi}\big(1+{\bar\varphi}|d_k|\big)\chiup_1 - \big(\sum_{j=2}^D{\hat A}_{1j}{\hat\psi}'_j\big)r^2 + {\bar\varphi}{\bar\psi}\chiup_1^2 \Bigg)Y_k,\\
			\xi_{k+1}= \xi_k + r^2,\qquad t_{k+1}= t_k + r^2.
			\end{array}\right.
			\end{eqnarray}}
		
		\ELSE{}
		
		\IF{$|d_k|\leq r\sqrt{\lambda_1({\bf X}_k,T-t_k)}$}
		
		\STATE{Let $r_{k+1}= |d_k|/\parallel\sigma_k^T{\bf N}_k\parallel$ (tangent ellipsoid to hyperplane, i.e. ${\Sigma_k\approx \widetilde{\Sigma}_k}$)}
		
		\ELSE{}
		
		\STATE{Let $r_{k+1}=r$ (regular ellipsoid $E_k^r$)}
		\ENDIF
		
		\STATE{Let ${\bm \omegaup_{k+1}}\sim{\cal S}_D$, ${\bm \mu}=\sigma^{-1}{\bf b}$ and take one step inside $\Omega$ \cite[algorithm 3.7]{Milstein_Tretyakov_Book}:}
		\begin{equation}
		\label{F:One_step_WoE}
		\left\{
		\begin{array}{l}
		{\bf X}_{k+1}= {\bf X}_k + \sigma_k{\bm \omegaup_{k+1}}r_{k+1}, \\
		Y_{k+1}= Y_k + Y_kc_k r_{k+1}^2/D +
		Y_k{\bm \mu}^T{\bm \omegaup_{k+1}}r_{k+1}, \\
		Z_{k+1}= Z_k + Y_kf_k r_{k+1}^2/D,\\ 
		\xi_{k+1}= \xi_k,\,\, t_{k+1}= t_k+r^2_{k+1}/D.
		\end{array}
		\right.
		\end{equation}
		
		\ENDIF
		
		\STATE{Update $d_{k+1}$ and $\big({\bf N}_{k+1},{\bf X}_{k+1}^{\partial\Omega}\big)$ according to \eqref{eq:N_and_proj}, and let $k=k+1$}
		
		\STATE{Safeguard: {\bf If} $d_{k+1}>0$ (overshoot) {\bf then} let ${\bf X}_{k+1}={\bf X}_{k+1}^{\partial\Omega}$ and $d_{k+1}=0$, {\bf end}}
		
		\ENDWHILE 
	\end{algorithmic}
	\label{A:WoE}
\end{algorithm}

\section{Detailed implementation}\label{S:Implementation}
\setlength\parindent{0pt}

The proofs of all lemmas in this Section are compiled in Appendix \ref{Ap:Proofs}.

\subsection[Largest eigenvalue]{Gershgorin estimation of the largest eigenvalue}\label{SS:Gershgorin}
In order to avoid overshooting the boundary with a hop from ${\bf X}_k$, an upper bound of $\lambda_1({\bf X}_k,T-t)$ is needed.
	In the event that there is no closed formula for it, accurately determining $\lambda_1({\bf X}_k,T-t_k)$ (for every timestep of every trajectory in the Monte Carlo simulation) may add up to a prohibitive computational overhead. For that reason, Milstein's WoE takes a uniform upper bound 	
	$\tilde{\lambda}_1\geq\max_{{\bf x}\in\Omega}\max_{0\leq t<T}\lambda_1({\bf x},T-t)$. However, $\tilde{\lambda}_1$ itself may not be straightforward to estimate, either; or unnecessarily large, thus requiring too small an $h$ in order to reach the required accuracy. In Algorithm \ref{A:WoE}, we propose substituting it by the non-uniform, sharper upper bound given below, with no detriment to $\delta$.

\begin{lemma}
	\label{Th:Gershgorin}
	$\lambda_1({\bf x},T-t)$ can be bounded above at ${\cal O}(D^2)$ cost by
	\begin{equation}
	\label{eq:Gershgorin}
	\lambda_1({\bf x},T-t)\leq \max_{j=1,\ldots,D}\sum_{i=1}^D|a_{ij}({\bf x},T-t)|.
	\end{equation}
\end{lemma}

The cost of this approach is thus ${\cal O}(D)$ cheaper than extracting the spectrum.

\subsection{Drift removal}\label{SS:Girsanov}
In order to remove the drift from (\ref{eq:Milstein_sys}), ${\bm\mu}=\sigma^{-1}{\bf b}$ must be calculated, which may be costly, specially in high dimension---unless $\sigma$ is lower triangular, where ${\bm\mu}$ can be computed easily by forward substitution at cost ${\cal O}(D^2)$. 
In the solution of BVPs, $\sigma$ can always be chosen lower triangular as the Cholesky matrix of $A({\bf x},t)$, because $A$ is positive definite. 

On the other hand, imagine that the data is a non-lower-triangular diffusion matrix ${\hat\sigma}$, and the point of the Monte Carlo calculation is to compute population densities such as those mentioned in Section \ref{S:Representation} (of the mean first exit time, etc.). Then, ${\hat\sigma}$ should be replaced in (\ref{eq:Milstein_sys}) by the (lower triangular) Cholesky matrix of 
${\hat\sigma}{\hat\sigma}^T)$.  In that scenario, it may be critical to obtain the lower triangular $\sigma$ analytically before start, for computing the Cholesky factorization numerically at every time step involves an ${\cal O}(D^3)$ cost per step \cite{Gill_matrix_updates}. In the remainder of the paper, we shall assume without loss of generality that $\sigma$ is lower triangular.

\subsection[Distance map]{Construction of the distance map}\label{SS:Eikonal}

Determining the triple of boundary data $\{d({\bf x}),{\bf x}^{\partial\Omega},{\bf N}({\bf x})\}$ for a point ${\bf x}$ may be time-consuming when $\Omega$ has a nontrivial shape. In general, the distance function (or distance map) inside $\Omega$ obeys the Eikonal equation in ${\mathbb R}^D$
\begin{equation}
\label{eq:Eikonal}
\parallel\nabla d({\bf x})\parallel= 1,\qquad u({\bf x}\in\partial\Omega)=0, \qquad d({\bf x}\in\Omega/\partial\Omega)<0.
\end{equation}

With the above convention that distances are negative inside $\Omega$, it holds
\begin{equation}
\label{eq:N_and_proj}
{\bf N}({\bf x})= \nabla d({\bf x}),\qquad {\bf x}^{\partial\Omega}= {\bf x}-d({\bf x})\nabla d({\bf x}).
\end{equation}  

Like in all schemes for bounded SDEs, the distance map must be calculated numerically prior to the Monte Carlo simulation. As proposed in \cite{Bernal2014}, Sethian's Fast Marching is the method of choice.

\subsection[Tangent ellipsoid]{Construction of the tangent ellipsoid}\label{SS:Ellipsoid}
At a given location ${\bf X}_k$ sufficiently far from the boundary, WoE draws ${\bf X}_{k+1}$ from the surface of the ellipsoid $E^r_k$, which is inscribed in the ball $B\big({\bf X}_k,r\lambda_1^{1/2}({\bf X}_k,T-t_k)\big)$.
When $|d_k|<r\lambda_1^{1/2}({\bf X}_k,T-t_k)$, there is no guarantee that $E^r_k$ does not stick out of $\partial\Omega$.
Milstein's WoE postulates that 
the tangent ellipsoid $\Sigma_k$ be taken  in that case (see \cite[algorithm 3.7]{Milstein_Tretyakov_Book}).
This prevents boundary overshoots while maximising the probability of ${\bf X}_{k+1}$ being absorbed, thus leading to the least average number of hops---namely ${\cal O}(1/h)$, see \cite[section 6.4.3]{Milstein_Tretyakov_Book}---and hence to an optimally efficient algorithm. (Even though $\Sigma_k$ might be larger than $E_r({\bf X}_k)$, the asymptotic weak convergence rate is ${\cal O}(h)$, as proven in \cite{Milstein_WoE}.)

Therefore, the determination of $\Sigma({\bf x})$ is needed for implementing Milstein's WoE. This is now formally solved by Lemma \ref{Th:Anisotropic} below---seemingly, a new result. 

\begin{lemma}
	\label{Th:Anisotropic}
	Let $\Omega$ be a closed domain in ${\mathbb R}^D$ (not necessarily smooth).
	Then, $\Sigma_k=E_k^{\Psi({\bf X}_k)}=E_k^{\overline{\rho}}$, where $\Psi({\bf x})$ is the solution of the anisotropic Eikonal equation 
	\begin{equation}
	\label{eq:AnisotropicA}
	||\nabla{\Psi}^TA_k\nabla{\Psi}||= 1\textrm{ and }\Psi>0 \textrm{ in }{\overline{\Omega}},\qquad \Psi=0 \textrm{ on }\partial\Omega.
	\end{equation}
\end{lemma}

Lemma \ref{Th:Anisotropic} clarifies a fundamental issue of WoE. When $A$ is a constant matrix, \eqref{eq:AnisotropicA} need be solved just once before the simulation, and then $\Psi({\bf x})$ will be evaluated in computing time---analogously to the Eikonal equation for the distance map. In fact, the solution to \eqref{eq:AnisotropicA} is required only inside a narrow shell on the inner side of $\partial\Omega$. A Fast-Marching-like method for Lemma \ref{Th:Anisotropic} can solve \eqref{eq:AnisotropicA} only there without regard to the rest of $\Omega$, thus cutting back on preprocessing overhead.

There are, however, two caveats to using Lemma \ref{Th:Anisotropic}. Numerical methods for the anisotropic Eikonal equation are
less developed (this is further commented on in Section \ref{S:Conclusions}).  
Moreover, in the event of a non-constant matrix $A({\bf x},T-t)$, solving \eqref{eq:AnisotropicA} at every ${\bf X}_k$ which needs it will in general be out of the question. 

For those reasons, we introduce the straightforward approximation of Lemma \ref{Th:Tangent_ellipsoid}, which in Section \ref{S:Experiments} is shown to work very well. 
The idea is to replace $\Sigma_k$ by the ellipsoid tangent to the closest tangent hyperplane, which we call ${\widetilde\Sigma}_k$.

\begin{lemma}
\label{Th:Tangent_ellipsoid}
The ellipsoid centred at ${\bf X}_k$ and tangent to the hyperplane tangent to $\partial\Omega$ 
at ${\bf X}_k^{\partial\Omega}$ 
is given by ${\widetilde \Sigma}_k=\{{\bf y}\,|\,({\bf y}-{\bf X}_k)^TA^{-1}_k({\bf y}-{\bf X}_k)=\rho^2_{k+1}\}$, where 
\begin{equation}
\label{eq:Half-space}
\rho_{k+1}=\frac{|d_k|}{||\sigma_k^T{\bf N}_k||}.
\end{equation}
When $\Big(\partial\Omega\cap B({\bf X}_k,|d_k|)\Big)\subset \Pi_k$,  the ellipsoid $E^{r_{k+1}}_k$ is fully inside $\Omega$, where
\begin{equation}
\label{eq:Rule_r}
r_{k+1}=\min{(\rho_{k+1},r)}.
\end{equation}
\end{lemma}

 Lemma \ref{Th:Tangent_ellipsoid} yields the required value $r_{k+1}$ for the next hop (check lines 18-22 in Algorithm \ref{A:WoE}). (Note that $A_k^{-1}$ is never used.) The cost of this construction is ${\cal O}(D^2)$, due to the product $\sigma_k^T{\bf N}_k$.

When the tangency point between $\Pi_k$ and $\widetilde{\Sigma}_k$ lies outside of $\Omega$, there is a nonzero probability that ${\bf X}_{k+1}$ overshoots: this is the reason of the safeguard in Algorithm \ref{A:WoE}, line 26. If $\partial\Omega$ is smooth at ${\bf X}_k^{\partial\Omega}$, the portion of $\widetilde{\Sigma}_k$ sticking out of $\partial\Omega$ tends to zero as $h\to 0^+$ (i.e. as $r\to 0^+$), since $\Pi_k\to \partial\Omega$ around ${\bf X}_k^{\partial\Omega}$. However, this may not be the case close to cusps or corners, specially if $A_k$ is a very eccentric ellipsoid with the major semiaxis parallel to $\Pi_k$. In that worst case scenario, it is important that $\widetilde{\Sigma}_k$ shrinks with $h$ so that the overshooting probability tends to vanish. This is the point of the rule $r_{k+1}=\min{(\rho_{k+1},r)}$ in Lemma \ref{Th:Tangent_ellipsoid}: Algorithm \ref{A:WoE} hops on $\widetilde{\Sigma}_k$ only if it is smaller than $E^r_k$ (i.e. if the tangency point is nearby); otherwise it sticks to $E^r_k$, knowing (by Lemma \ref{Th:Tangent_ellipsoid}) that the probability of ${\bf X}_{k+1}$ overshooting goes asymptotically to zero. 
We note, however, that this is an heuristic reasoning rather than a rigourous proof that the construction given by Lemma \ref{Th:Tangent_ellipsoid} preserves $\delta=1$ from MM in presence of curved boundaries.

\subsection[Reflecting boundary]{Change of coordinates close to the reflecting boundary} \label{SS:Rotation}

\paragraph*{Notation.} In this subsection, column vectors in ${\mathbb R}^D$ and ${\mathbb R}^{D-1}$ are respectively written in bold (like ${\bf N}_k$) and with arrows (like ${\vec N}$). 
Given a $D\times D$ matrix (like $\sigma$), a starred matrix (like $\sigma_*$) denotes the submatrix obtained by removing the first column and row from the former.
 Matrix (or vector) elements superfluous for the discussion are depicted by $*$ (like ${\bf N}_k^T=(*,{\vec N}^T)$). $I_S$ stands for the identity matrix in dimension $S$; and ${\bf e}_n$ is the $n^{th}$ column of $I_D$. We drop subindex $k$ except for ${\bf X}_k$, ${\bf X}_k^{\partial\Omega}$, $d_k$, $\nabla \psi_k$, and ${\bf N}_k$. All numbers are real.\newline  

When the trajectory ${\bf X}_k$ is closer than $r$ to $\partial\Omega_R$, coordinates are locally changed so that the new origin is ${\bf X}_k$ and the new first component points towards $-{\bf N}_k$ (i.e. inwards). The rotation is thus defined as 
\begin{equation}
\label{F:Qdef_1}
Q({\bf X}_k^{\partial\Omega}-{\bf X}_k)= (-|d_k|,0,\ldots,0)^T,
\end{equation}
where $Q$ is an orthogonal matrix, i.e. $Q^T=Q^{-1}$ (not unique, in general). In order to perform the rotations in an efficient way, we adapt the approach in \cite{Buchmann_Petersen_SIAM}. 

The Givens matrix\footnote{Householder transformations could be used as well, see \cite{Buchmann_Petersen_SIAM}.} $G^{\theta}(i,j)$ (where $i<j$) is defined element-wise as
\begin{equation}
\label{eq:Givens_1st}
[G^{\theta}(i,j)]_{kl} = \delta_{kl}, \textrm{ except: } 
\left\{\begin{array}{ll}
[G^{\theta}(i,j)]_{ii}= \cos{\theta}, & [G^{\theta}(i,j)]_{ij}= \sin{\theta}, \\
\textrm{$[G^{\theta}(i,j)]_{ji}= -\sin{\theta}$}, & [G^{\theta}(i,j)]_{jj}=  \cos{\theta}.\\
\end{array}\right.
\end{equation}
Givens matrices are orthogonal and have two properties of interest to us:
\begin{enumerate}
	\item Let ${\bf v}=(v_1,\ldots,v_D)$ be an arbitrary vector. There is an angle $\theta({\bf v})$, given by
	
\begin{equation}
\label{eq:theta}
\left\{
\begin{array}{lll}
\cos{\theta({\bf v})}=1,&
\sin\theta({\bf v})=0&\textrm{ if $v_i=v_j=0$,}\\
\cos\theta({\bf v})= \frac{v_i}{\sqrt{v_i^2+v_j^2}},&
\sin\theta({\bf v})=\frac{v_j}{\sqrt{v_i^2+v_j^2}},
&\textrm{ otherwise,}
\end{array}
\right.
\end{equation}
	
	such that the action of $G^{\theta({\bf v})}(i,j)$ on a vector zeroes its $j^{th}$ element, may change the $i^{th}$ one, and leaves the rest unchanged:
	\begin{equation*}
	{\bf v}'=G^{\theta({\bf v})}(i,j){\bf v}=(v_1,\ldots,v_{i-1},v'_i,v_{i+1},\ldots,v_{j-1},v_j'=0,v_{j+1},\ldots,v_D)^T.
	\end{equation*}
	\item If $v_i=v_j=0$, and ${\bf v}'=G^{\beta}(i,j){\bf v}$, then $v_i'=v_j'=0$ for any $\beta\in{\mathbb R}$.
\end{enumerate}

Let ${\bf v}^{(1)}$ be a vector and ${\bf v}^{(2)}=G^{\theta_2}(1,2){\bf v}^{(1)}$, where $\theta_2=\theta({\bf v}^{(1)})$ according to (\ref{eq:theta}). Then, by construction, the second element of ${\bf v}^{(2)}$ is zero and every other one is unchanged except for the first one. Similarly, the second and third elements of ${\bf v}^{(3)}=G^{\theta_{3}}(1,3){\bf v}^{(2)}$ with $\theta_3=\theta({\bf v}^{(3)})$ per (\ref{eq:theta}) are zero. Iterating, it is clear that
\begin{equation}
\label{eq:Givens_order} 
G^{\theta_D}(1,D)\cdots
G^{\theta_{3}}(1,3)
G^{\theta_{2}}(1,2){\bf v}^{(1)}=:G{\bf v}^{(1)}= (\beta,0,\ldots,0)^T.
\end{equation} 
Performing the rotations in the above order has the following properties:

\begin{lemma}
\label{Th:Properties_of_Q}
Let $||{\bf v}^{(1)}||>0$, and $G$ be the matrix product of the Givens rotations $G=\Pi_{k=D}^2G^{\theta_k}(1,k)$ as in (\ref{eq:Givens_order}), i.e. zeroing the vector elements from the second to the last, starting from ${\bf v}^{(1)}$. Then: (i) $G_*$ is lower triangular; (ii) $\beta=||{\bf v}^{(1)}||>0$.
\end{lemma}

In order for $Q{\bf N}_k=(-1,0,\ldots,0)$, we set ${\bf v}^{(1)}={\bf N}_k$, sequentially construct $G^{\theta_2}(1,2),\ldots,G^{\theta_D}(1,D)$, and set
\begin{eqnarray}
\label{eq:Q_of_Givens}
Q=  -G^{\theta_{D}}(1,D)\cdots G^{\theta_{2}}(1,2),
\end{eqnarray}
and therefore
\begin{align}
\label{eq:Q^T_of_Givens}
Q^{-1}&= -[G^{\theta_{2}}(1,2)]^{-1}\cdots 
[G^{\theta_{D}}(1,D)]^{-1}
\nonumber\\
&=
-[G^{\theta_{2}}(1,2)]^T\cdots 
[G^{\theta_{D}}(1,D)]^T
\nonumber\\
&=-[G^{-\theta_{2}}(1,2)]\cdots 
[G^{-\theta_{D}}(1,D)].
\end{align}

In order to premultiply a vector, $Q$ (or $Q^T$) need not be formed, but the two sets $\{\cos{\theta_2},\ldots,\cos{\theta_D}\}$ and $\{\sin{\theta_2},\ldots,\sin{\theta_D}\}$ are calculated and stored in advance, and later used whenever needed. To signify that the sequence of Givens rotations is performed on ${\bf v}$, without ever constructing $Q$, we shall write
\begin{equation}
\label{F:Q_matrix_vector}
Q{\bf v}= -\verb|ROTATIONS|({\bf v},\theta_2\shortrightarrow\theta_D),\qquad
Q^T{\bf v}= -\verb|ROTATIONS|({\bf v},-\theta_D\shortrightarrow -\theta_2).
\end{equation}

The cost of carrying out the $D-1$ Givens rotations sequentially is just ${\cal O}(D-1)$, instead of ${\cal O}\big((D-1)^2\big)$ as would be the case for an explicit matrix-vector multiplication. Several vectors in Algorithm \ref{A:WoE} need to be rotated while computing the reflection. For instance,
$({\hat\psi}'_1,\ldots,{\hat\psi}'_D)^T=Q\nabla\psi_k=$\verb|-ROTATIONS|$(\nabla\psi_k,{\theta_2\shortrightarrow\theta_D)}$.

Let us now focus on the following two parts of matrix $A$ in the rotated frame:
\begin{equation}
\label{eq:A_gorro}
{\hat A}=QAQ^T=
\left(\begin{array}{cc}
{\hat A}_{11}=(QAQ^T)_{11} & * \\ * &{\hat A}_*= (QAQ^T)_*
\end{array}\right).
\end{equation} 
 Because $A$ is positive definite, so are $QAQ^T$ (since $Q$ is orthogonal) and\footnote{If $M$ is positive definite, $M_*$ too, since $0<{\bf y}^TM{\bf y}={\vec y}^TM_*{\vec y}$, for any ${\bf y}=(0,{\vec y}^T)^T$ such that $||{\vec y}||>0.$} $(QAQ^T)_*$, and $(QAQ^T)_{11}>0$. 
Then ${\hat A}_{11}>0$ and ${\hat A}_{11}r^2+d_k^2>0$ in Algorithm \ref{A:WoE}. 

\paragraph*{Updating the decomposition of the rotated submatrix.} Let us now address another important computational aspect, not discussed in MM. It is required to compute the matrix-vector product ${\Lambda_*}{\vec\nu}$, where  
$\Lambda_*$ is defined by
\begin{equation}
{\hat A}_*= \Lambda_*\Lambda_*^T.
\end{equation}

Thanks to the positive-definiteness of ${\hat A}_*$, $\Lambda_*$ could be obtained by Cholesky factorization at a cost ${\cal O}\big((D-1)^3\big)$. This is taxing if ${\bf N}_k$ is not constant on $\partial\Omega_R$ and/or in high dimension. Fortunately, the factorization $A=\sigma\sigma^T$ can be efficiently recycled into $\Lambda_*$---although, as it will be shown next, not necessarily in the standard way.
The following notation will be convenient:

\begin{equation}
\label{F:Split_of_Q_and_sigma}
Q=
\left(\begin{array}{cc}
Q_{11} & {\vec q}_1^T \\ {\vec q}_2 & Q_*
\end{array}\right),
\qquad
\sigma=
\left(\begin{array}{cc}
\sigma_{11} & {\vec 0}^T \\ {\vec s} & {\sigma}_*
\end{array}\right).
\end{equation} 

Since $\sigma$ is nonsingular and lower triangular by construction (see Section \ref{SS:Girsanov}), $\sigma_{11}\neq 0$, $\det{\sigma_*}\neq 0$, and $\sigma_*$ is lower triangular. On the other hand, $Q_*$ need not be orthogonal or even nonsingular. Writing out ${\hat A}= Q\sigma\sigma^TQ^T$ gives

\begin{equation}
\label{F:Update}
\Lambda_*\Lambda_*^T= Q_*\sigma_*\sigma_*^TQ_*^T + {\vec w}{\vec w}^T, \qquad \textrm{where }{\vec w}= \sigma_{11}{\vec q}_2 + Q_*{\vec s}.
\end{equation}

Therefore, $\Lambda_*\Lambda_*^T$ can be regarded as the factorization $(Q_*\sigma_*)(Q_*\sigma_*)^T$ plus the rank-one update ${\vec w}{\vec w}^T$. Borrowing a standard ansatz
from \cite{Gill_matrix_updates}, we set

\begin{equation}
\label{F:Cholesky_ansatz}
\Lambda_* = Q_*\sigma_*(I_{D-1}+\alpha{\vec z}{\vec z}^T)
\end{equation}  

and look for suitable $\alpha$ and ${\vec z}$. Inserting (\ref{F:Cholesky_ansatz}) into (\ref{F:Update}) and noting that $(I_{D-1}+\alpha{\vec z}{\vec z}^T)(I_{D-1}+\alpha{\vec z}{\vec z}^T)^T=$ 
$(I_{D-1}+\alpha{\vec z}{\vec z}^T)^2=I_{D-1}+{\vec z}(\alpha^2{\vec z}^T{\vec z} + 2\alpha){\vec z}^T$, one has

\begin{equation}
(\alpha^2{\vec z}{\vec z}^T + 2\alpha)(Q_*\sigma_*{\vec z})(Q_*\sigma_*{\vec z})^T= {\vec w}{\vec w}^T,
\end{equation}

which can be readily solved by letting
$\alpha=\alpha_+$, where
\begin{equation}
\label{eq:alpha}
\alpha^2{\vec z}^T{\vec z} + 2\alpha=1 \Rightarrow \alpha_{\pm}= 
\frac{-1\pm\sqrt{1+||{\vec z}||^2}}{||{\vec z}||^2},
\end{equation}

(we choose $\alpha_+$ for concreteness). Then, the vector ${\vec z}$ is given by

\begin{equation}
\label{eq:Update_solution}
Q_*\sigma_*{\vec z}= {\vec w}.
\end{equation}

Formulas \eqref{F:Cholesky_ansatz}-\eqref{eq:Update_solution} allow for calculation of $\Lambda_*\nu$ at ${\cal O}\big(D(D-1)\big)$ cost. However, Lemma \ref{Th:Update_old} below shows that this ansatz may fail.

\begin{lemma}
\label{Th:Update_old}
Let ${\bf N}_k^T=(N_1,{\vec N}^T)$. The update formula $Q_*\sigma_*{\vec z}=\sigma_{11}{\vec q}_2 + Q_*{\vec s}$ (with $Q_*$, ${\vec q}_2$, $\sigma_{11}$, $\sigma_*$, and ${\vec s}$ from \eqref{F:Split_of_Q_and_sigma}) is inconsistent if and only if $N_1=0$. In that case, $\det{Q_*}=0$. Otherwise, there is one unique ${\vec z}$, given by
\begin{equation}
\sigma_*{\vec z}= {\vec s}-\frac{\sigma_{11}}{N_1}{\vec N}.
\end{equation}
\end{lemma}  

We stress that Lemma \ref{Th:Update_old} is independent of the way in which the rotation is implemented. Swapping the first component for the first nozero one of ${\bf N}_k$ in \eqref{F:Qdef_1} does not help, either: the permutation matrices involved induce structural changes (the analogous of $\sigma_*$ is no longer lower triangular, for instance), with the result that Milstein's formulas in Algorithm \ref{A:WoE} would have to be reworked. 

When the first component of ${\bf N}_k$ is zero, $\Lambda_*$ can still be computed at ${\cal O}\big((D-1)^2\big)$ cost by  using Algorithm \verb|ZCHUD| in LINPACK \cite{Dongarra_LINPACK} (implemented as \verb|cholupdate| in Matlab).
This exploits Lemma \ref{Th:Properties_of_Q}: since $Q_*\sigma_*$ is lower triangular, $\Lambda_*$ can be seen as a rank-one update of a Cholesky factorization. Despite the fact that $Q_*\sigma_*\sigma_*^TQ_*^T$ is only positive semidefinite if $N_1=0$, \verb|ZCHUD| would also work in that case---thanks to the fact that one Cholesky matrix, namely $Q_*\sigma_*$, is available in the first place. However, using \verb|ZCHUD| would involve forming the full matrices. \newline

Instead, we put forward the following analytical approach, specifically tailored to the case $N_1=0$. It relies on the following decomposition:

\begin{equation}
\label{eq:rectangular_decomp}
\Lambda_*\Lambda_*^T=
\Big[Q_*\sigma_*\,\big|\,{\vec w}\Big]
\left[\begin{array}{c}
\sigma_*^TQ_*^T \\ \smallskip \\ {\vec w}^T
\end{array}\right]=
\Big[Q_*\sigma_*\,\big|\,{\vec w}\Big]
ZZ^T
\left[\begin{array}{c}
\sigma_*^TQ_*^T \\ \smallskip \\ {\vec w}^T
\end{array}\right],
\end{equation}

where $Z$ is an orthogonal matrix of order $D$. If it can be chosen such that

\begin{equation}
\label{eq:Z_split}
Z^T
\left[\begin{array}{c}
\sigma_*^TQ_*^T \\ \smallskip \\ {\vec w}^T
\end{array}\right]=
\left[\begin{array}{c}
H \\ \smallskip \\ {\vec 0}^T
\end{array}\right],
\end{equation}

then $\Lambda_*=H^T$. It turns out that $Z$ can be found analytically, leading to a closed formula for $\Lambda_*$.

\begin{lemma}
\label{Th:Z}
Let ${\bf N}_k^T=(0,{\vec N}^T)$; $Q_*$, ${\vec q}_2$, $\sigma_{11}$, $\sigma_*$, and ${\vec s}$ from \eqref{F:Split_of_Q_and_sigma},
and ${\vec w}$ from \eqref{F:Update}. Then,
\begin{equation}
\label{eq:Z_ap}
Z=
\left[\begin{array}{cc}
I_{D-1}-\frac{1}{||{\vec r}||^2}{\vec r}{\vec r}^T & \frac{{\vec r}}{||{\vec r}||} \\
&\\
\frac{{\vec r}^T}{||{\vec r}||} & 0
\end{array}\right],
\end{equation}
where $\sigma_*{\vec r}={\vec N}$. Furthermore,  $\Lambda_*=Q_*\sigma_*+\frac{1}{||{\vec r}||}{\vec w}{\vec r}^T$.
\end{lemma}

\paragraph*{Remark.}  $\Lambda_*$ in Lemma \ref{Th:Z} is not necessarily triangular---but note that this was not required, anyways (check \cite[formula 6.19]{Milstein_Tretyakov_Book}, where $\Lambda_*$ is called $\lambda$).\newline

Let ${\vec v}$ be an arbitrary column vector in ${\mathbb R}^{D-1}$. Observe that
\begin{eqnarray}
\label{eq:Qv}
Q\left(\begin{array}{c}
0 \\ {\vec v}
\end{array}\right)=
\left(\begin{array}{c}	
* \\ Q_*{\vec v}
\end{array}\right)= &-\verb|ROTATIONS|\Big( 
\left(\begin{array}{c}	
0 \\ {\vec v}
\end{array}\right)
,\theta_2\rightarrow \theta_D \Big),\\
Q^T\left(\begin{array}{c}
0 \\ {\vec v}
\end{array}\right)=
\left(\begin{array}{c}	
* \\ Q_*^T{\vec v}
\end{array}\right)= &-\verb|ROTATIONS|\Big( 
\left(\begin{array}{c}	
0 \\ {\vec v}
\end{array}\right)
,-\theta_D\rightarrow -\theta_2 \Big).
\end{eqnarray}

Let $[QM]_n$ be the $n^{th}$ column of $QM$, where $M$ is a square matrix of order $n$. Then, $[QM]_n$ can be calculated as $[QM]_n=(QM){\bf e}_n=-\verb|ROTATIONS|(M{\bf e}_n,\theta_2\rightarrow \theta_D)$. Furthermore, the first column of  matrix ${\hat A}$ is $QAQ^T{\bf e}_1=-Q(\sigma\sigma^T{\bf N}_k)$, i.e.
\begin{equation}
\label{F:A1}
\left(\begin{array}{c}
{\hat A}_{11} \\ \vdots \\ {\hat A}_{D1}
\end{array}\right)=
{\hat A}{\bf e}_1
= 
\verb|ROTATIONS|\Big(
\sigma\sigma^T{\bf N}_k
,\theta_2\shortrightarrow\theta_D\Big).
\end{equation}

The concrete calculation of $\Lambda_*\nu$ with cost ${\cal O}\big(D(D-1)\big)$ is listed as Algorithm \ref{A:Update}. Complexity is dominated by the product $\sigma_*{\vec\nu}$ and the forward substitutions.

\begin{algorithm}[h!]
	\caption{Fast computation of $\Lambda_*{\vec{\nu}}$ in Algorithm \ref{A:WoE} (line 13)}
	\begin{algorithmic}[1]
		\STATE{{\bf Data:}  $\sigma$, ${\vec\nu}$, Givens angles $\{\cos{\theta_2,\ldots,\sin{\theta_D}}\}$, ${\bf N}_k=(N_1,{\vec N})^T$, tiny $tol>0$}
		\STATE{Let $\left(\begin{array}{c}
			* \\ {\vec q}_2 
			\end{array}\right)=
			-\verb|ROTATIONS|\Big(
			(1,0,\ldots,0)^T
			,\theta_2\rightarrow\theta_D
			\Big)$}
		\STATE{Compute $Q_*(\sigma_*{\vec\nu})$  and $Q_*{\vec s}$ by \eqref{eq:Qv} and let ${\vec w}=\sigma_{11}{\vec q}_2+Q_*{\vec s}$}
		\IF{$|N_1|>tol$}
\STATE{Compute ${\vec z}$ from $\sigma_*{\vec z}={\vec s}-\frac{\sigma_{11}}{N_1}{\vec N}$ by forward substitution}
\STATE{ $\Lambda_*{\vec\nu}=Q_*(\sigma_*{\vec\nu})+\big(\sqrt{1+||{\vec z}||^2}-1\big)\,\frac{{\vec z}^T{\vec \nu}}{||{\vec z}||^2}\,{\vec w}$}		
		\ELSE
		\STATE{Compute ${\vec r}$ from $\sigma_*{\vec r}={\vec N}$ by forward substitution}
		\STATE{$\Lambda_*{\vec\nu}=
		Q_*(\sigma_*{\vec{\nu}})+	
\frac{{\vec r}^T{\vec\nu}}{||{\vec r}||}{\vec w}$}
		\ENDIF		
	\end{algorithmic}
	\label{A:Update}
\end{algorithm}

\section{Numerical experiments}\label{S:Experiments}

In this section, we report numerical results obtained with Algorithm \ref{A:WoE}---henceforth, they are labelled as "MM''. The algorithm of Constantini {\em et al.} has also been used for comparison (``REF''). (Specifically, the version with pseudonormal variables, which has a weak convergence rate ${\cal O}(h^{1/2-\epsilon}),\epsilon>0$ \cite{Constantini1998}.) Both codes have been written in Matlab, and run on a laptop. 
In that case (even though Algorithm \ref{A:WoE} is written in sequential form), the codes should be fully vectorized, which is critical for speed (check appendix \ref{Ap:Code}). 

The Matlab code and data files used for this paper are available at the journal repository, and upon request from the author. 

Let
${\phi}_h^{(1)},\ldots,{\phi}_h^{(N)}$ be $N$ iid numerical approximations to $\phi$ in  \eqref{eq:Parabolic_representation}. The numerical approximation to the Feynman-Kac functional is 
\begin{equation}
u_{h,N}({\bf x}_0,t)= \frac{1}{N}\sum_{j=1}^N \phi_h^{(j)}.
\end{equation}
Asymptotically (i.e. as $N\to\infty$ and $h\to 0^+$),  $\varepsilon_{N,h}\sim {\cal N}(Ch^{\delta},{\mathbb V}[\phi_h]/N)$, where $\varepsilon_{N,h}$ is the root mean-square (RMS) error of $\phi_h$, $C>0$ is a constant and $\delta>0$ is the weak convergence rate of the scheme \cite{Milstein_Tretyakov_Book}. In all the ensuing experiments, the pointwise exact solution $u_{ex}({\bf x}_0,T)={\mathbb E}[\phi]\neq 0$ is known. Then, after setting an accuracy goal $\varepsilon$, $N$ is chosen accordingly, i.e. such that $2\sqrt{{\mathbb V}[\phi_h]/ N}= .20\times \varepsilon$. The idea is that asymptotically the RMS error carries less than a $20\%$ statistical error with a large probability. 
However, for better assessment and comparison across experiments, the relative error $\epsilon_h^{rel}:=|1-u_{h,N}({\bf x}_0,T)/u_{ex}({\bf x}_0,T)|$ is reported instead. 

\paragraph{Remark.} In practice, one starts simulating trajectories $j=1,2,...$, replaces ${\mathbb V}[\phi_h]$ by the sample estimate so far, $V_j$, and stops as soon as $\sqrt{V_j/j}\leq \varepsilon/10$.

\subsection[Example I]{ Example I (three-dimensional)}\label{SS:Example1}

We take this from \cite{Gobet&Menozzi_2010}. The coefficients of \eqref{eq:ParabolicBVP_mixedBCs} are:
$c=0$, ${\bf b}=(y,z,x)^T$, $\varphi=-||{\bf x}||^2$, 
\begin{eqnarray}
\label{F:ParabolicGM}
\sigma(x,y,z)=\left(\begin{array}{ccc}
\sqrt{1+|z|} & 0 & 0 \\
\frac{1}{2}\sqrt{1+|x|} & \sqrt{\frac{3}{4}} \sqrt{1+|x|} & 0 \\
0 & \frac{1}{2} \sqrt{1+|y|} &   \sqrt{\frac{3}{4}}\sqrt{1+|y|}		
\end{array}
\right).
\end{eqnarray}

By choosing as exact solution $u_{ex}({\bf x},t)=xyz$ (independent of time), and $g=p=u_{ex}$, the remaining coefficients $\psi$ and $f$ are derived from \eqref{eq:ParabolicBVP_mixedBCs}, i.e. $\psi={\bf N}^T\nabla u_{ex}-\varphi u_{ex}$, and 
$f= y^2z + z^2x + x^2y + \frac{1}{2}\sqrt{1+\mid z \mid}\sqrt{1+\mid x \mid}+\frac{\sqrt{3}}{2}x\sqrt{1+\mid x \mid}\sqrt{1+\mid y \mid}$.\newline

As in \cite{Gobet&Menozzi_2010}, we set $\Omega=B({\bf 0},1)$  and ${\bf x}_0=(.56,.52,.30)^T$; thus 
$u_{ex}({\bf x}_0,T)\approx.08736$. We consider three different sets of BCs: purely absorbing, purely reflecting, and mixed. In the latter case, the hemisphere with $z<0$ is absorbing and the other one, reflecting. (Note that the reflection is not conormal.) The convergence of 
$\epsilon_h^{rel}$ w.r.t. $h$ is shown in Table \ref{T:GM_ball}, along with the results with REF. Regardless of the BCs, MM is both more accurate and has a faster convergence rate (estimated by least-squares regression), in fact very close to the theoretical value $\delta=1$.

\begin{table}[h!]
	\centering	
	\begin{tabular}{|l|cc|cc|cc|}
\multicolumn{1}{c}{}&
\multicolumn{2}{c}{absorbing BCs ($T=\infty$)}&
\multicolumn{2}{c}{mixed BCs ($T=\infty$)}&
\multicolumn{2}{c}{reflecting BCs ($T=1$)}\\
\hline
$h$ & MM  & REF & MM  & REF & MM  & REF  \\
		\hline
.0128 &.07594&.33949&.29351&.27006&.26478&.21171 \\
.0064 &.04022&.25511&.15755&.20827&.13380&.16075 \\
.0032 &.02101&.18721&.08050&.15630&.07619&.12005 \\
.0016 &.00930&.13048&.04348&.11840&.03377&.09089 \\
.0008 &.00490&.10203&.01737&.08428&.01424&.06866 \\
.0004 &.00232&.07380&.01016&.06603&.00839&.05175 \\
.0002 &.00140&.04984&.00405&.04455&.00327&.04189 \\
.0001 &.00070&.04219&.00223&.03169& .00213&.03023 \\
		\hline
\multicolumn{1}{|c|}{$\delta$} &0.98&0.44&1.02&0.44&1.03&0.40\\
\hline
	\end{tabular}
	\caption{Convergence of $\epsilon_h^{rel}$ for Example I in $\Omega=B({\bf 0},1)$ at ${\bf x}_0=(.56,.52,.30)^T$. (Entries carry about $\pm 20\%$ statistical error.)}
	\label{T:GM_ball}
\end{table}

In many applications boundaries are not so smooth as on a sphere; we also solved this problem in the box $\Omega=[-\sqrt{2}/2,\sqrt{2}/2]^3$. We only show results for the purely absorbing and purely reflecting cases; see Table \ref{T:GM_in_box}. Computational times have been included.

Since the boundary is nonsmooth now (due to the corners), accuracy and convergence rate are bound to deteriorate. In particular, theoretical rates no longer apply. Nonetheless, MM is still the more efficient integrator. For instance---acording to Table \ref{T:GM_in_box}---REF took $340$ s. to attain a relative error of $.12722$, while MM took just $65.7$ s. for $.12573$ (with the same confidence interval). We stress that this is the case where MM performs worst (purely reflecting nonsmooth boundary).

\begin{table}[h!]
	\centering	
	\begin{tabular}{|l|cl|cl|cl|cl|}
		\multicolumn{1}{c}{}&
		\multicolumn{4}{c}{absorbing BCs ($T=\infty$)}&
		\multicolumn{4}{c}{reflecting BCs ($T=1$)}\\
		\hline
$h$ & MM  & time & REF & time & MM  & time & REF & time\\
		\hline
.0128 & .03400& 	0.85&	.08591&	0.67&		.45085&		5.78&		.70420&	5.34\\
.0064 & .01605&		7.34&	.06711&	1.14&	.27996&		14.3&		.53835&	12.6\\
.0032 & .00770&	50.7&	.05659&	1.97&		.18981&	28.9&		.40311&	21.2\\
.0016 & .00488&	279&	.04218&	3.66&	.12573&	65.7&		.30115&	40.1\\
.0008 & .00221& 2957&	.03130&	13.7&		.10147&	115&	.23036&	88.0\\
.0004 & .00106&  25715&	.02441&	24.2&	.07350&	245&	.15882&	178\\
.0002 & .00062&	145236&	.01930&	72.6&	.05324&	987&	.12722&	340\\
.0001 &.00027&849117& .01086&	520&	.04752&	3076&	.08534&	644\\
		\hline
		\multicolumn{1}{|c|}{$\delta$}&0.98&&0.40&&			0.47&&0.43&	\\
		\hline
	\end{tabular}
	\caption{Convergence of $\epsilon_h^{rel}$ (within about $\pm 20\%$ statistical error) for Example I in $\Omega=[-\sqrt{2}/2,\sqrt{2}/2]^3$ at ${\bf x}_0=(.56,.52,.30)$. Times in s. on a single processor.}
\label{T:GM_in_box}
\end{table}

\subsection[Example II]{Example II (arbitrary-dimensional)}\label{SS:Example2}

This problem features crossed second derivatives, oscillating coefficients and solution, and can be posed in any dimension. The diffusion matrix is
\begin{equation}
\label{F:Tril}
\sigma=\left(\begin{array}{cccc}
1 & 0  & \ldots & 0 \\
1 & 1  & \ldots & 0 \\
\vdots & \vdots & \ddots & \vdots \\
1 & 1  & \ldots & 1
\end{array}\right)
\Rightarrow
\sigma\sigma^T= A= 
\left(\begin{array}{ccccc}
1 & 1 & 1 & \ldots & 1 \\
1 & 2 & 2 & \ldots & 2 \\
1 & 2 & 3 & \ldots & 3 \\
\vdots & \vdots & \vdots & \ddots & \vdots \\
1 & 2 & 3 & \ldots & D 
\end{array}\right).
\end{equation}
The exact solution is $u_{ex}=\cos{\sum_{i=1}^D x_i}$. The coefficients are: ${\bf b}=(\sin{\pi x_1},\ldots,\sin{\pi x_D})^T$, $c=0$, $\varphi=-1$, and $g=u_{uex}$ ($f$ and $\psi$ can be derived from \eqref{eq:ParabolicBVP_mixedBCs}). We take $T=\infty$ (an elliptic PDE, so that no initial condition is needed), and $\Omega=B({\bf 0},1)$ in ${\mathbb R}^D$.\newline

\begin{table}[h!]
	\centering	
	\begin{tabular}{|l|cc|c|cc|}	
		\hline	
$h$ & MM (D=4) & MM (D=5) & REF (D=5) & MM (D=6)& MM (D=9)  \\
		\hline
.0064& .00829&.02507&.19833&.07132&.45630\\
.0032&.00358& .01320&.14175&.04110&.20576\\
.0016& .00132& .00655&.10550&.01727&.11527\\
.0008&.00067& .00371&.08002&.00914&.06331\\
.0004&.00032& .00162&.04708&.00512&.02679\\
.0002& ---		& .00087&.03566&.00207&.01468\\
.0001& ---		&.00043&.02916&.00131&.00642\\
.00005&---&---&.01832&.00061&.00399\\ 
.000025&---&---&.01248&---&.00169\\
		\hline
$\delta$&1.18&0.98&0.49&0.99&1.00\\ 
		\hline
	\end{tabular}
\caption{Convergence of $\epsilon_h^R$ (within about $\pm20\%$ statistical error) for the Example II in increasing dimension $D$. Solved in $\Omega=B({\bf 0},1)$ at $T=\infty$ (i.e. it is an elliptic problem) and at ${\bf x}_0=(-\pi/4,0,\ldots,0)$. BCs of mixed type (absorbing if $z<0$). Missing entries were not calculated because they would take unacceptably long.}
\label{T:Radcliffe_in_ball}
\end{table}

\begin{figure}[htb]
	\centering
	\includegraphics[width=1\linewidth]{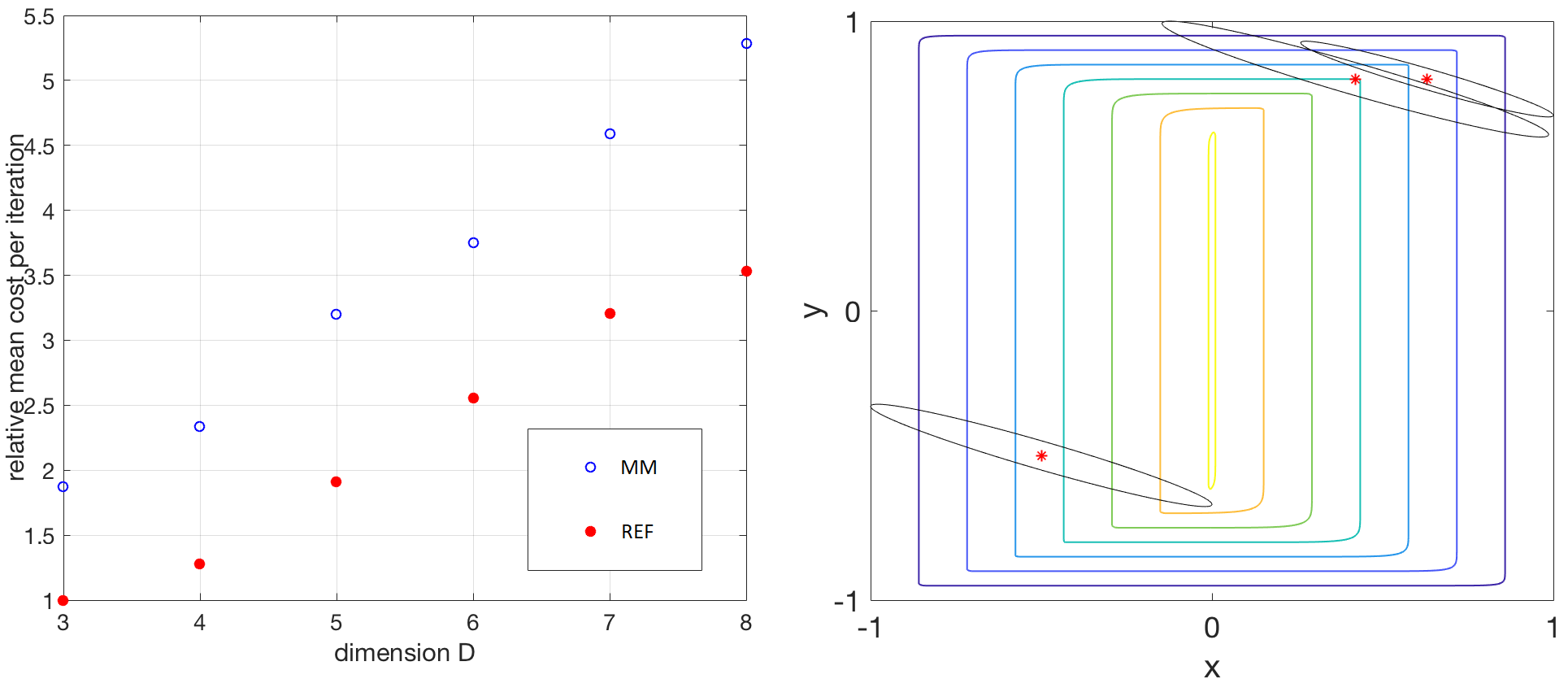}
\caption{(Left) Relative average cost per iteration of MM and REF. (Right) Contours of the solution to \eqref{eq:Eiko_sardina}. Three tangent ellipses are overlaid, for which the major semiaxis is proportional to $\Psi({\bf x})$.}
\label{F:Cost_iter}
\end{figure}

In Table \ref{T:Radcliffe_in_ball}, we consider the same mixed BCs as before (i.e. reflecting on the upper hemisphere). The solution with MM at ${\bf x}_0=(-0.7854,0,\ldots,0)^T$ is evaluated for increasing $D$. For comparison, the solution with REF in $D=5$ is also included. Times are not included, but let us give an example: REF took $3119$ s. for a relative accuracy of $.01248$, while MM (always in $D=5$) took $50.9$ s. for $\epsilon_h^{rel}=.01320$. As usual, the difference grows with increasing accuracy (as ${\cal O}(1/\epsilon_h^{rel})$, see \cite{Giles_y_yo}) due to the better $\delta$ of MM. 

We use Example II to compare the cost per iteration of MM and REF. Specifically, we solve it in B({\bf 0},1) with $T=1$, same ${\bf x}_0$, $N=10^5$ (fixed),  $h=.001$ and purely reflecting BCs (so that the number of steps is roughly constant). The total times taken by MM and REF are shown on the left of  Figure \ref{F:Cost_iter}, for $D=3,\ldots,8$. This indicates that the mean cost per iteration of MM is less than twice as with REF.

All numerical experiments so far have made use of Gershgorin's circles theorem (Section \ref{SS:Gershgorin}) for estimating $\lambda_1$. We checked that doing this does not have a discernible impact on accuracy (compared with using the exact $\lambda_1$), while it cuts back on computational overhead (by about $15\%-20\%$, in these examples).

\begin{figure}[htb]
\centering
\includegraphics[width=1.1\linewidth]{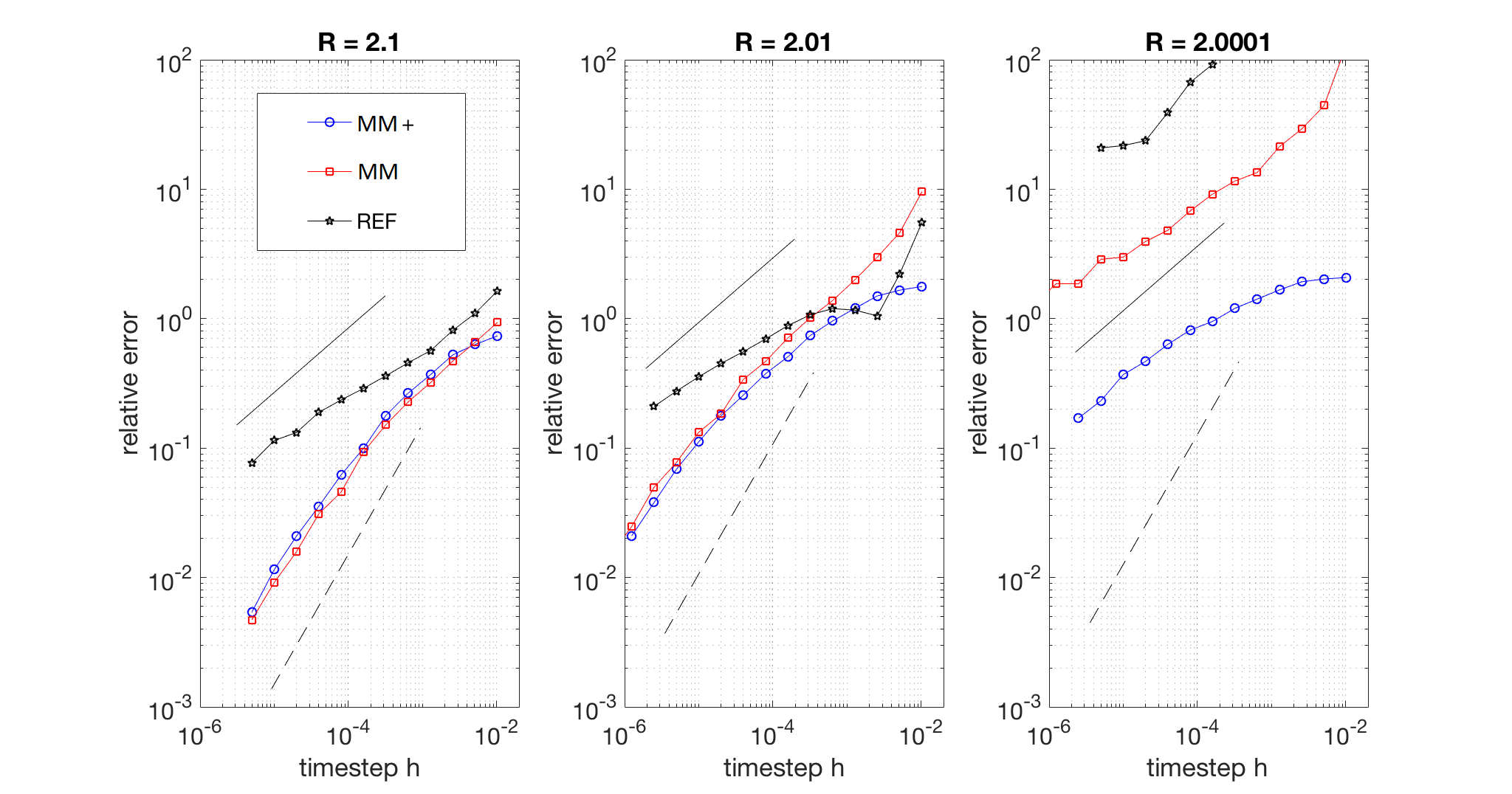}
\caption{Convergence of $\epsilon_h^{rel}$ (within $\pm 20\%$ statistical error) for Problem III and increasing difficulty at the corners. The dashed and solid segments (with slope $1$ and $1/2$, respectively) are meant as reference. With $R=2.0001$ neither MM nor MM+ have yet attained the asymptotic linear regime.}
\label{F:Sardina}
\end{figure}

\subsection[Example III]{Example III (two-dimensional with sensitive corners)}\label{SS:Example3} 
The purpose of this problem is to investigate the connection between MM and the anisotropic Eikonal equation established in this paper. We pick $T=\infty,\,{\bf x}_0=(.823,-.875)^T$ and purely absorbing BCs in $\Omega=[-1,1]^2$. The coefficients are:
\begin{equation}
\label{eq:matrix_exp3}
A=
\left(\begin{array}{cc}
8 & -2.71 \\ -2.71 & 1\\
\end{array}\right)
=\sigma\sigma^T
\qquad\Rightarrow\qquad
\sigma\approx
\left(\begin{array}{cc}
2.8284 & 0 \\ -0.9581 & 0.2863\\
\end{array}\right),
\end{equation}
$c=0$, ${\bf b}=(\sin{\pi x},\sin{\pi y})^T$, $g=u_{ex}$, and $f$ is derived from the exact solution,
\begin{equation}
\label{eq:Sardina}
u_{ex}(x,y)=\frac{1}{x^2+y^2-R}.
\end{equation}
The point of \eqref{eq:Sardina} is that, as $R\to 2$, the exact solution tends to infinity on the four corners, substantially compounding the difficulty of integrating the corresponding SDEs in the domain with nonsmooth boundary. 
Problem III is solved with MM and REF, as before, plus a new version of MM where, instead of using Lemma \ref{Th:Tangent_ellipsoid} for $\widetilde{\Sigma}_k$, $\Sigma_k$ is calculated exactly according to Lemma \ref{Th:Anisotropic}---we call this integrator ``MM+''. In order to do so, the anisotropic Eikonal equation  
\begin{equation}
\label{eq:Eiko_sardina}
||\nabla{\Psi}^TA\nabla{\Psi}||= 1\textrm{ and }\Psi>0 \textrm{ in }{\overline{\Omega}},\qquad \Psi=0 \textrm{ on }\partial\Omega,
\end{equation}
is solved before the Monte Carlo simulation. (Importantly, note that thanks to the fact that $A$ is constant, only one such equation must be solved.) We do so using the Fast Sweeping scheme described in \cite{Fast_Sweeping}, and store the numerical solution in a $251\times 251$ grid, from where $\Psi({\bf X}_k)$ is later interpolated ``on the fly''.  (Lines 19-21 in Algorithm \ref{A:WoE} must then be replaced by $r_{k+1}=\Psi({\bf X}_k)$.)

The contours of $\Psi$ are sketched on Figure \ref{F:Cost_iter} (right). Figure \ref{F:Sardina} shows the convergence of $\epsilon_h^{rel}$ for increasing difficulty: $R=\{2.1,2.01,2.0001\}$. Note that MM+ converges faster than MM, and both versions can tackle this problem much better than REF---which in the hardest case $R=2.0001$,  essentially cannot cope.

\section{Conclusions}\label{S:Conclusions}

Numerical experiments 
show a very satisfactory performance of Algorithm \ref{A:WoE}, in terms of accuracy, speed, and weak convergence rate.  
We stress that it has a favourable complexity, rendering it suitable for high-dimensional problems.

The connection, established in this paper, with the anisotropic Eikonal equation, paves the way for the efficient numerical integration of general absorbed diffusions in domains with corners. (The only other method with this capability known to the author is the Walk on Rectangles algorithm \cite{Deaconu2006}, which is nonetheless restricted to the Brownian motion and some variations thereof with constant coefficients.) On the other hand, solving the anisotropic Eikonal equation is nontrivial and the object of intense current research. Tailoring existing solvers for the anisotropic Eikonal equation to the particular features of the algorithm presented here is left as future work.

Finally, we point out that the most promising direction for further development of this algorithm is producing a Multilevel formulation, along the lines of \cite{Giles_y_yo}.

\section*{Acknowledgements}
The author thanks J. A. Acebr\'on for proposing this work and  helpful discussions throughout. Portuguese FCT funding is gratefully acknowledged. 

\appendix
\section{Proofs of Lemmas in Section 4}\label{Ap:Proofs}

{\bf Proof of Lemma \ref{Th:Gershgorin}.} This is a direct application of Gershgorin's circles theorem: for any eigenvalue $\lambda$ of any $D\times D$ matrix $M$, there is $1\leq i \leq D$ such that
\begin{equation}
|\lambda-M_{ii}|\leq \sum_{i\neq j}|M_{ij}|.
\end{equation}
Since $A$ is positive definite, $\lambda>0$. By the triangle inequality, $|\lambda-a_{ii}|\geq |\lambda|-|-a_{ii}|=\lambda-|a_{ii}|$, and hence
\begin{equation}
\label{F:Gershgorin}
\lambda_1({\bf x},T-t)\leq \max_{i=1,\ldots,D}\sum_{j=1}^D|a_{ij}({\bf x},T-t)|.
\end{equation}
The cost is that of adding up $D$ elements for each of the $D$ rows, and comparing those $D$ values to find the largest one. \hfill$\Square$\newline

{\bf Proof of Lemma \ref{Th:Anisotropic}}. Let $M$ be a $D\times D$ symmetric positive definite matrix, with spectrum $M{\bf m}_i=\mu_i{\bf m}_i$,  where $\mu_1\geq \ldots \geq \mu_D>0$ and $||{\bf m}_i||=1$. The (positive) viscosity solution of the anisotropic Eikonal equation
\begin{equation}
\label{eq:AnisotropicM}
||\nabla{\Psi}^TM\nabla{\Psi}||= 1 \textrm{ in }{\overline{\Omega}},\qquad \Psi=0 \textrm{ on }\partial\Omega,
\end{equation}
at ${\bf x}\in\Omega$, namely $\Psi({\bf x})$, can be interpreted as the arrival time at ${\bf x}$ of a monotonic front which marches anisotropically with velocity $\mu_i$ along the direction ${\bf m}_i$ (see \cite{Schwenke2012} and references therein). 

By symmetry, it is also the arrival time on $\partial\Omega$ of a front marching outwards from ${\bf x}$ with the same velocities in the opposite directions (i.e. $\mu_i$ along $-{\bf m}_i$).

Then, at time $t'>0$, the front stemming from ${\bf x}$ is the ellipsoid with semiaxes of length $\mu_1t',\ldots,\mu_Dt'$ along the orthogonal directions ${\bf m}_1,\ldots,{\bf m}_D$. Since $M$ and $M^{-1}$ have the same eigenvectors
and inverse eigenvalues, that ellipsoid is the locus of $({\bf y}-{\bf x})^TM^{-1}({\bf y}-{\bf x})=t'^2$. 

In particular, at time $t''=\Psi({\bf x})$, the ellipsoidal front arrives on the boundary for the first time, meaning that it is tangent to it. Therefore,
\begin{equation}
\Sigma_M({\bf x})= \{{\bf y}\,|\,({\bf y}-{\bf x})^TM^{-1}({\bf y}-{\bf x})=\Psi^2({\bf x})\}.
\end{equation} 
Setting $M=A_k$ and $\Sigma_k=E^{|\Psi({\bf X}_k)|}_{A({\bf X}_k,T-t_k)}({\bf X}_k)$ yields the desired tangent ellipsoid.\hfill$\Square$\newline

{\bf Proof of Lemma \ref{Th:Tangent_ellipsoid}.} Without loss of generality, let us take the origin at ${\bf X}_k$. Set $m({\bf y},\rho)= {\bf y}^TA^{-1}_k{\bf y}-\rho^2$.  
The sought-for ellipsoid is then $\widetilde{\Sigma}_k=\{{\bf y}\,|\,m({\bf y},\rho_{k+1})=0\}$. The tangency point ${\bf y}_0$ belongs both to $\widetilde{\Sigma}_k$ and to $\Pi_k$, so that for some $q\neq 0$,
\begin{equation}
\label{eq:elipf}
\nabla m({\bf y}_0,\rho_{k+1})=2A^{-1}_k{\bf y}_0=q {\bf N}_k.
\end{equation}
Since $\det A_k\neq 0$, \eqref{eq:elipf} implies
\begin{equation}
{\bf y}_0=\frac{q}{2}A_k{\bf N}_k\Rightarrow
\left\{
\begin{array}{l}
\frac{q^2}{4}{\bf N}_k^TA_k^{T}A^{-1}_kA_k{\bf N}_k= \rho^2_{k+1},
\textrm{ and }\\
{\bf N}_k^T{\bf y}_0=\frac{q}{2}{\bf N}_k^TA_k{\bf N}_k.
\end{array}
\right.
\end{equation}
Moreover, $A_k=A^T_k$ and hence
\begin{equation}
\rho_{k+1}= \frac{{\bf N}_k^T{\bf y}_0}{\sqrt{{\bf N}_k^TA^T_k{\bf N}_k}}=
\frac{|d_k|}{||\sigma^T{\bf N}_k||}.
\end{equation}

The condition $\Big(\partial\Omega\cap B({\bf X}_k,|d_k|)\Big)\subset \Pi_k$ is equivalent to $\partial\Omega$ being planar inside $B({\bf X}_k,|d_k|)$. Let us assume it is. Then, $E^r_k$ sticks out of $\partial\Omega$ iff it sticks out of $\Pi_k$. If it doesn't, then $\rho_{k+1}\geq r$. If it does, then $\widetilde{\Sigma}_k=E^{\rho_{k+1}}_k$ is concentric to and inside of $E^r_k$, so that $\rho_{k+1}<r$. In both cases, $E^{r_{k+1}}({\bf X}_k)$ given by rule \eqref{eq:Rule_r} is inside $\Omega$. \hfill$\Square$\newline

{\bf Proof of Lemma \ref{Th:Properties_of_Q}.} In preparation, let $\theta\in{\mathbb R}$, $1<j,k\leq D$, and ${\bf w}=G^{\theta}(1,k){\bf e}_j$ . Then, 
\begin{equation}
\textrm{if $k\neq j$ and $j>1$} \Rightarrow G^{\theta}(1,k){\bf e}_j= {\bf e}_j.
\end{equation}
To see this, simply note that $w_1=({\bf e}_j)_1\cos{\theta}+({\bf e}_j)_k\sin{\theta}$ and $w_k=-({\bf e}_j)_1\sin{\theta}+({\bf e}_j)_k\cos{\theta}$ are both zero if $j>1$ and $j\neq k$, while the other elements of ${\bf w}$ are unaffected by the rotation and thus are the same as in ${\bf e}_j$.

Assume now $1<i<j\leq D$. Then,
\begin{align}
\label{eq:Gij}
G_{ij}&= {\bf e}_i^TG{\bf e}_j \nonumber\\
	     &= {\bf e}_i^T G^{\theta_D}(1,D)\cdots G^{\theta_{i+1}}(1,i+1) G^{\theta_i}(1,i) \cdots G^{\theta_2}(1,2) {\bf e}_j \nonumber\\
&=\Big([G^{\theta_{i+1}}(1,i+1)]^T\cdots[G^{\theta_D}(1,D)]^T{\bf e}_i\Big)^T \Big(G^{\theta_i}(1,i) \cdots G^{\theta_2}(1,2) {\bf e}_j\Big)\nonumber\\
&=\Big(G^{-\theta_{i+1}}(1,i+1)\cdots G^{-\theta_D}(1,D){\bf e}_i\Big)^T{\bf e}_j\nonumber\\
&={\bf e}_i^T{\bf e}_j=0.		     
\end{align}

In \eqref{eq:Gij}, we have used the facts that: $[G^{\theta}(m,n)]^T=G^{-\theta}(m,n)$; that $\Pi_{k=i}^2G^{\theta_k}(1,k){\bf e}_j={\bf e}_j$ because the interval $2,\ldots,i$ does not include $j$, by assumption; and analogously $i+1,\ldots,D$ does not include $i$. 
Since $G_{1<i\leq D,1<j<i}=0$, $G_*$ is lower triangular, thus proving (i).

On the other hand, the upper row of $G$ is not zero, in general. In fact, working out those entries (which is tedious but straightforward), one has:
\begin{align}
\label{eq:G1n}
G_{11}&= \cos{\theta_2}\cos{\theta_3}\cos{\theta_4}\cdots\cos{\theta_D},\nonumber\\
G_{1k}&=\sin{\theta_k}\cos{\theta_{k+1}}\cos{\theta_{k+2}}\cdots\cos{\theta_D},\,(2\leq k\leq D).
\end{align}

Let us now prepare for (ii). First, let ${\bf v}$ be a vector and ${\bf v}'=G^{\theta({\bf v})}(1,k){\bf v}$. By \eqref{eq:theta},
\begin{equation}
v'_1=\sqrt{\big(v_1^{(1)}\big)^2+\big(v_k^{(k)}\big)^2}\geq 0.
\end{equation}

Consequently, in the sequence ${\bf v}^{(k)}=G^{\theta_k}(1,k){\bf v}^{(k-1)}$ with $k=2,\ldots,D$, it holds that $v_1^{(2)}\geq 0,\ldots,v_1^{(D)}\geq 0$. By \eqref{eq:theta},  $\cos{\theta_3},\ldots,\cos{\theta_D}$ are thus all nonnegative. (Only $\cos{\theta_2}<0$ iff $v_1^{(1)}<0$). 

Clearly, $\beta^2=({\bf v}^{(1)})^TG^TG{\bf v}^{(1)}=
({\bf v}^{(1)})^T{\bf v}^{(1)}= ||{\bf v}^{(1)}||^2>0$, by hypothesis. 
Also, it is always possible to pick $n$ such that $v^{(1)}_n\neq 0$. Then, $v_n^{(1)}={\bf e}_n^T{\bf v}^{(1)}=
{\bf e}_n^TG^TG{\bf v}^{(1)}=
{\bf e}_n^TG^T(\beta,0,\ldots,0)^T=
\beta G_{1n}$. This prevents $G_{1n}$ from being zero, so that $\beta=v_n/G_{1n}$. Moreover, by the previous discussion and \eqref{eq:G1n}, $G_{1n}\neq 0$ implies that $\cos{\theta_{n+1}}\cdots\cos{\theta_D}>0$. 

Therefore, it is clear that, if $n=1$, sign($G_{11}$)=sign($\cos{\theta_2}$)=sign($v_1^{(1)}$), by \eqref{eq:theta} and \eqref{eq:G1n}.

If $1<n\leq D$, sign($G_{1n}$)=sign($\sin{\theta_n}$)=sign($v_n^{(n-1)}$), by \eqref{eq:theta}. Let us show that $v_n^{(n-1)}=v_n^{(1)}$. This occurs because ${\bf v}^{(k)}=G^{\theta_k}(1,k){\bf v}^{(k-1)}$ preserves all but the first and $k^{th}$ elements of ${\bf v}^{(k-1)}$, so that if $n>k$, then $v_n^{(k)}=v_n^{(k-1)}$. Iterating, it is clear that ${\bf v}^{(n-1)}=G^{\theta_{n-1}}(1,n-1)\cdots G^{\theta_{2}}(1,2){\bf v}^{(1)}$ has preserved the elements $v^{(n-1)}_n,\ldots,v^{(n-1)}_D$ from ${\bf v}^{(1)}$. Then, by \eqref{eq:theta},  sign($\sin{\theta_n}$)=sign($v_n^{(n-1)}$)=sign($v_n^{(1)}$).

Summing up, we have proved that sign($\beta$)= sign($v^{(1)}_n/G_{1n})>0$. In fact, $\beta=+||{\bf v}^{(1)}||$, which is (ii). \hfill$\Square$\newline

{\bf Proof of Lemma \ref{Th:Update_old}.}  We recall that $\det{\sigma_*}\neq 0$ and $\sigma_{11}\neq 0$ because $\sigma\sigma^T$ is positive definite. By definition $Q(N_1,{\vec N}^T)^T=(-1,{\vec 0}^T)^T$ and hence
\begin{equation}
\label{eq:Q*N_up}
Q_{11}N_1+{\vec q}_1^T{\vec N}=-1,
\end{equation}
\begin{equation}
\label{eq:Q*N_down}
N_1{\vec q}_2+Q_*{\vec N}={\vec 0}.
\end{equation}
Writing out $QQ^T=I_D$ and $Q^TQ=I_D$ in terms of the blocks defined in  (\ref{F:Split_of_Q_and_sigma}) yields
\begin{equation}
\label{eq:QQ1}
Q_{11}^2 + {\vec q}_1^T{\vec q}_1= 1,
\end{equation}
\begin{equation}
\label{eq:QQ2}
Q_{11}^2 + {\vec q}_2^T{\vec q}_2= 1,
\end{equation}
\begin{equation}
\label{eq:QQ3}
-Q_{11}{\vec q}_2 = Q_*{\vec q}_1,
\end{equation}
\begin{equation}
\label{eq:QQ4}
-Q_{11}{\vec q}_1 = Q_*^T{\vec q}_2,
\end{equation}
\begin{equation}
\label{eq:QQ5}
{\vec q}_2{\vec q}_2^T + Q_*Q_*^T= I_{D-1},
\end{equation}
\begin{equation}
\label{eq:QQ6}
{\vec q}_1{\vec q}_1^T + Q_*^TQ_*= I_{D-1}.
\end{equation}

\begin{itemize}
	
	\item Let us start by the 'if' part. Assume that $N_1=0$. 
	
	Since $||{\bf N}_k||=1$, then ${\vec N}_k\neq {\vec 0}$. From (\ref{eq:Q*N_down}), $Q_*{\vec N}={\vec 0}$, and therefore $\det{Q_*}=0$. Applying Sylvester's determinant identity to (\ref{eq:QQ5}) gives:
	\begin{equation}
	\label{eq:Sylvester}
	\det{(Q_*Q_*^T)}= \det{(I_{D-1}-{\vec q}_2{\vec q}_2^T)}= 1-{\vec q}_2^T{\vec q}_2 \,\Rightarrow\, {\vec q}_2^T{\vec q}_2=1.
	\end{equation}
	
	Then, by \eqref{eq:QQ2}, $Q_{11}=0$. By \eqref{eq:QQ3}, this means that ${\vec q}_2^TQ_*={\vec 0}^T$, so that ${\vec q}_2$ is not in the range of $Q_*$ (i.e. it is orthogonal to the column space of $Q_*^T$). On the other hand, $Q_*(\sigma_*{\vec z})$ has a projection on ${\vec q}_2$, since ${\vec q}_2^TQ_*(\sigma_*{\vec z})=\sigma_{11}\neq 0$. Consequently, there cannot be such $\sigma_*{\vec z}$ if $N_1=0$.
	
	\item Assume now that the system is inconsistent. Since $\sigma_{11}\neq 0$ and $Q_*{\vec s}$ in ${\vec w}$ clearly belongs to the range of $Q_*$, it follows that ${\vec q}_2\neq 0$ does not. This implies ${\vec q}_2^TQ_*={\vec 0}^T$, and by \eqref{eq:Q*N_down}, $N_1=0$.  This proves the ``only if'' part.
	
	\item Finally, assume there is a solution. It has been proved that then $N_1\neq 0$ and $Q_{11}\neq 0$. By Sylvester, $\det{(Q_*Q_*^T)}=1-{\vec q}_2^T{\vec q}_2=Q_{11}^2$, whence $\det{Q_*}\neq 0$.  This means that $\sigma_*{\vec z}$, and thus ${\vec z}$, is unique.
	
	If $N_1\neq 0$, then $-\frac{\sigma_{11}}{N_1}Q_{\vec N}=\sigma_{11}{\vec q}_2$ by \eqref{eq:Q*N_down}. It follows that $Q_*\big( {\vec s}-\frac{\sigma_{11}}{N_1}{\vec N} \big)={\vec w}$. \hfill$\Square$
	\end{itemize}

{\bf Proof of Lemma \ref{Th:Z}.} Let us try to fit \eqref{eq:Z_split} with

\begin{equation}
\label{eq:Z_candidate}
\left[\begin{array}{cc}
\Phi & {\vec r} \\ {\vec r}^T & 0
\end{array}\right]
\left[\begin{array}{c}
\sigma_* ^TQ_*^T \\ \smallskip \\ {\vec w}^T
\end{array}\right]
=
\left[\begin{array}{c}
\Phi\sigma_*^T Q_*^T + {\vec r}{\vec w}^T\\
{\vec r}^T\sigma_*^T Q_*^T
\end{array}\right].
\end{equation}

In order for the bottom row of the rightmost matrix above to vanish, $Q_*\sigma_*{\vec r}={\vec 0}$.  Since $\det{\sigma_*}\neq 0$, this has a nontrivial solution iff $\det{Q_*}=0$, i.e. iff $N_1=0$, by Lemma \ref{Th:Update_old}. In that case, by \eqref{eq:Q*N_down}, $\sigma_*{\vec r}=\lambda{\vec N}$, with $\lambda\in{\mathbb R}$. (For definiteness, we pick $\lambda=1$.)

The second requirement for a candidate $Z$ is that it be orthogonal. By symmetry of the righmost matrix in \eqref{eq:Z_candidate}, it suffices to check that $Z^TZ=I$:

\begin{equation}
Z^TZ=
\left[\begin{array}{cc}
\Phi\Phi^T+{\vec r}{\vec r}^T & \Phi{\vec r} \\ {\vec r}^T\Phi^T & {\vec r}^T{\vec r}
\end{array}\right]
=
\left[\begin{array}{cc}
I_{D-1} & {\vec 0} \\ {\vec 0}^T & 1
\end{array}\right]
\Leftrightarrow
\left\{
\begin{array}{ll}
\Phi\Phi^T=I_{D-1}-{\vec r}{\vec r}^T & (i) \\
\Phi{\vec r}={\vec 0} & (ii) \\
||{\vec r}||=1 & (iii)
\end{array}	
\right.	
\end{equation}

Assume for the time being that $||{\vec r}||=1$. Let us take $\Phi=I_{D-1}-{\vec r}{\vec r}^T$. Then
$\Phi=\Phi^T$ and $\Phi\Phi^T=\Phi^2=I_{D-1}-2{\vec r}{\vec r}^T+{\vec r}({\vec r}^T{\vec r}){\vec r}^T=I_{D-1}-{\vec r}{\vec r}^T=\Phi$, meeting (i). For (ii), note that $\Phi{\vec r}=(I_{D-1}-{\vec r}{\vec r}^T){\vec r}={\vec r}-({\vec r}^T{\vec r}){\vec r}={\vec 0}$. Finally, condition (iii) is met simply by substituting ${\vec r}/||{\vec r}||$ for ${\vec r}$ (this makes the result independent of $\lambda$.)

Therefore, $Z$ in \eqref{eq:Z_ap} is orthogonal, and yields $\Lambda_*=Q_*\sigma_*(I_{D-1}-\frac{1}{||{\vec r}||^2}{\vec r}{\vec r}^T)+\frac{1}{||{\vec r}||}{\vec w}{\vec r}^T$. Replacing ${\vec r}$ by $\sigma_*^{-1}{\vec N}$ and using the fact that $Q_*{\vec N}={\vec 0}$ finishes the proof. \hfill$\Square$

\section{Matlab code snippets}\label{Ap:Code}

Vectorized functions are given below for \verb|ROTATIONS| (Section \ref{SS:Rotation}), forward subtitution with a lower triangular system, and distribution ${\cal B}$ (Algorithm \ref{A:WoE}).

\begin{footnotesize}
\begin{verbatim}
function varargout= rotations(opcion,varargin)
% [C,S]= rotations('init',vector) meaning that Givens(C,S)*vector=[1,0,...0]
% Q*vector= -rotations('forth',vector,C,S)
% Q'*vector= -rotations('back',vector,C,S) 
X1= deal(varargin{1}); [N,dim]= size(X1);
switch opcion
case 'init'
C= zeros(N,dim-1); S= zeros(N,dim-1); %stored in same order as they are defined: ORD
for k=2:dim
wk= sqrt(X1(:,1).^2+X1(:,k).^2); C(:,k-1)= X1(:,1)./wk; S(:,k-1)= X1(:,k)./wk;
cero= find(wk==0); if ~isempty(cero), C(cero,k-1)= 1; S(cero,k-1)= 0; end
X2= X1; %and now overwrite two components:
X2(:,1)= C(:,k-1).*X1(:,1) + S(:,k-1).*X1(:,k);
X2(:,k)= -S(:,k-1).*X1(:,1) + C(:,k-1).*X1(:,k);
X1= X2; clear X2
end %vector X's only nonzero entry is the first one now. 
varargout= {C,S}; return
case 'forth', rotations= [2:1:dim]; s1= +1; s2= -1; %order ORD 
case 'back',  rotations=[dim:-1:2]; s1= -1; s2= +1; %order ORD inverse
end
%Implement Givens rotations:
[C,S]= deal(varargin{2:3});
for k=rotations 
X2= X1; %and now overwrite two components:
X2(:,1)=    C(:,k-1).*X1(:,1) + s1*S(:,k-1).*X1(:,k);
X2(:,k)= s2*S(:,k-1).*X1(:,1) +    C(:,k-1).*X1(:,k);
X1= X2;
end
varargout= {X2}; %rotated X1
\end{verbatim}

\begin{verbatim}
function x= forsuvec(t,X,b) %vectorized forward substitution L*x=b
global PDEfile
[N,dim]= size(X); x= NaN(N,dim-1); %result
for k=1:dim-1
	fila= feval(PDEfile,t,X,'sigmarow',k+1); fila= fila(:,2:end);
	if k>1, x(:,k)= ( b(:,k) - sum(fila(:,1:k-1).*x(:,1:k-1),2) )./fila(:,k);
	else x(:,k)= b(:,1)./fila(:,1); end
end
\end{verbatim}

\begin{verbatim}
function [Eta]= distribucionboluda(N,dim)
Eta= randn(N,dim); R= sqrt(sum(Eta.^2,2)); Eta= Eta./(R*ones(1,dim));
\end{verbatim}

\end{footnotesize}

\end{document}